\documentclass[12pt,amstex]{amsart}

\usepackage{mathptmx}
\usepackage{mathrsfs}

\usepackage{verbatim}
\usepackage{url}
\usepackage[all]{xy}
\usepackage{color}

\usepackage[colorlinks=true,citecolor=blue, backref]{hyperref}

\usepackage{stmaryrd}
\usepackage{epsfig}
\usepackage{amsmath}
\usepackage{amssymb}

\usepackage{amscd}
\usepackage{graphicx}
\usepackage{pstricks}
\usepackage{tocvsec2}

\newtheorem{ques}{Question}

\newtheorem{thm}{Theorem}[section]
\newtheorem{lem}[thm]{Lemma}

\newtheorem{prop}[thm]{Proposition}
\newtheorem{cor}[thm]{Corollary}

\newtheorem{thma}{Theorem}

\theoremstyle{definition}
\newtheorem{de}[thm]{Definition}
\newtheorem{exam}[thm]{Example}
\theoremstyle{remark}
\newtheorem{rem}[thm]{Remark}

\numberwithin{equation}{section}

\def \N {\mathbb{N}}

\def \Z {\mathbb{Z}}
\def \R {\mathbb{R}}
\def \Q {\mathbb{Q}}

\def \O {\mathcal{O}}

\def \B {\mathcal B}

\def \X {\mathcal{X}}
\def \Y {\mathcal{Y}}

\def \id {{\rm id}}

\def \a {\alpha }
\def \b {\beta}
\def \ep {\epsilon}
\def \d {\delta}
\def \D {\Delta}

\def \diam {\rm diam}
\def \ra {\rightarrow}
\def \T {\mathbb T}

\begin{document}

\title{On systems disjoint from all minimal systems}

\author[W.~Huang]{Wen Huang}
\address[W. Huang]{School of Mathematical Sciences, University of Science and Technology of China, Hefei, Anhui 230026, China}
\email{wenh@mail.ustc.edu.cn}

\author[S.~Shao]{Song Shao}
\address[S. Shao]{School of Mathematical Sciences, University of Science and Technology of China, Hefei, Anhui 230026, China}
\email{songshao@ustc.edu.cn}

\author[H.~Xu]{Hui Xu}

\address[H. Xu]{Department of Mathematics, Shanghai Normal University, Shanghai, 200234, China}
\email{huixu@shnu.edu.cn}

\author[X.~Ye]{Xiangdong Ye}
\address[X. Ye]{School of Mathematical Sciences, University of Science and Technology of China, Hefei, Anhui 230026, China}
\email{yexd@ustc.edu.cn}

\thanks{This research is supported by National Key R$\&$D Program of China (No. 2024YFA1013601, 2024YFA1013602, 2024YFA1013600), and National Natural Science Foundation of China (12426201, 12371196, 12201599, 12031019 and 12090012).}

\begin{abstract}
Recently, G\'{o}rska, Lema\'{n}czyk, and  de la Rue characterized the class of automorphisms disjoint from all ergodic
automorphisms. Inspired by their work, we provide several characterizations of systems that are disjoint from all minimal systems.

For a topological dynamical system $(X,T)$, it is disjoint from all minimal systems if and only if there exist minimal subsets $(M_i)_{i\in\mathbb{N}}$ of $X$ whose union is dense in $X$ and each of them is disjoint from $X$
(we also provide a measure-theoretical analogy of the result).
For a semi-simple system $(X,T)$, it is disjoint from all minimal systems if and only if there exists a dense $G_{\delta}$ set $\Omega$ in $X \times X$ such that for every pair $(x_1,x_2) \in \Omega$, the subsystems $\overline{\O}(x_1,T)$ and $\overline{\O}(x_2,T)$ are disjoint. Furthermore, for a general system a characterization similar to the ergodic case is obtained.
\end{abstract}

\maketitle


\section{Introduction}\label{sec-intro}

In the ring of integers $\Z$ two integers $m$ and $n$ have no common factor if whenever $k|m$ and $k|n$ then $k = \pm 1$. They are disjoint if whenever $m|k$ and $n|k$, then also $mn|k$. In $\Z$ these two notions coincide. In his seminal
paper in 1967, Furstenberg \cite{Fur67} introduced the notion of disjointness in both ergodic theory and topological dynamics to compare distinct systems. In this paper, we provide several characterizations of systems that are disjoint from all minimal systems.

\subsection{Background}

Recently G\'{o}rska, Lema\'{n}czyk, and  de la Rue \cite{GLR24} describe the class Erg$^\perp$ of automorphisms disjoint from all ergodic
automorphisms. It is shown that a measure preserving system $(X,\B,\mu,T)$ is disjoint from all ergodic
automorphisms if and only if in the space $(\bar{X},\bar \mu)$ of ergodic components of $(X,\B,\mu,T)$, for $\bar{\mu}\times \bar{\mu}$-a.e. pair $(\bar{x}_1,\bar{x}_2)$ the ergodic
transformations $(T_{\bar{x}_1}, T_{\bar{x}_2})$ are disjoint. This result shows that Erg$^\perp$ is a much richer class
than what was previously thought. 
Soon after Glasner and Weiss \cite{GW24} provided an alternative proof of this result and their proof works just as well for any countable acting group. Inspired by their work, we provide several characterizations of topological systems that are disjoint from all minimal systems.

First we recall some basic notions. A {\em topological dynamical system}, or simply a {\em system}, is a pair $(X,T)$, where $X$ is a compact metric space with metric $\rho$ and $T: X \rightarrow X$ is a homeomorphism. We say a system $(X, T)$ is \emph{non-trivial} if $X$ is not a singleton. A system $(X,T)$ is called {\em minimal} if $X$ contains no proper non-empty closed invariant subsets (a point of $x\in X$ is called a {\it minimal point}, if the orbit closure of $x$ is minimal). A system $(X,T)$ is {\em transitive} if there is some point $x\in X$ with $\overline{\{T^nx:n\in\mathbb{Z} \}}=X$ (such a point is called a {\it transitive point}); it is {\em weakly mixing} if the product system $(X\times X,T\times T)$ is transitive.

For two systems $(X,T)$ and $(Y,S)$, a \emph{joining} of $(X,T)$ and $(Y,S)$ is a closed subset of the product $X \times Y$ that remains invariant under the joint action of $T \times S$ (i.e. $(T\times S)J= J$ \footnote{If one considers continuous maps instead of homeomorphisms, then it only requires $(T\times S) J\subset J$.}) and projects onto each coordinate space. If the only possible joining is the trivial product, that is, $X \times Y$ itself, then the systems $(X,T)$ and $(Y,S)$ are said to be {\em disjoint}, denoted by $(X,T)\perp (Y,S)$ or $X \perp Y$.

Note that if two systems are disjoint, then at least one of them must be minimal. In his seminal work \cite{Fur67}, Furstenberg showed that a weakly mixing system is disjoint from any minimal distal system.\footnote{A system $(X,T)$ is distal if $\inf_{n\in \Z}\rho(T^nx,T^ny)>0$ for all $x\neq y\in X$.} He also posed a problem to characterize the classes of systems that are disjoint from all minimal systems and distal systems (\cite[Problem G]{Fur67}). Petersen addressed this in \cite{Pet70}, concluding that the class of systems disjoint from all distal systems is comprised of weakly mixing minimal systems. Subsequently, lots of researchers have endeavored to characterize systems that are disjoint from all minimal systems.

Let $\mathcal{M}$ represent the class of minimal systems, and let $\mathcal{M}^{\perp}$ represent the class of systems that are disjoint from any minimal system. For system $(X,T)$, we also use $X\perp \mathcal{M}$ to denote that $(X,T)\in {\mathcal M}^\perp$, i.e. it is disjoint from all minimal systems.  Furstenberg showed in \cite{Fur67} that a totally transitive system\footnote{A system $(X,T)$ is {\em totally transitive} if $(X,T^n)$ is transitive for all $n\in \Z\setminus\{0\}$.} with dense periodic points is disjoint from all minimal systems. A direct corollary of this is that the Bernoulli shift is an element of $\mathcal{M}^{\perp}$. After that, many mathematicians have provided conditions under which a system belongs to $\mathcal{M}^{\perp}$. A necessary condition for transitive systems within $\mathcal{M}^{\perp}$ was first established by Huang and Ye in \cite{HY05}. They showed that a transitive system that is disjoint from all minimal systems must be weakly mixing and have dense minimal points.

The class $\mathcal{M}^{\perp}$ encompasses several types of systems, including:\footnote{For definitions of small periodic points and regular minimal points, refer to \cite{HY05}. As these concepts are not utilized in the subsequent discussion, their definitions are omitted here.}
\begin{itemize}
\item Totally transitive systems with dense sets of small periodic sets (\cite{HY05}),
\item Weakly mixing systems with dense sets of regular minimal points (\cite{HY05}),
\item Weakly mixing systems with dense sets of distal points (\cite{DSY12, Opr10}).
\end{itemize}
Moreover, the dynamical properties of a hyperspace system (see Subsection \ref{subsec-qf} for definitions) can also indicate disjointness. In \cite{LYY15}, the authors showed that if a system $(X,T)$ is weakly mixing and its hyperspace system $(2^{X}, T)$ has dense distal points, then $(X,T)$ belongs to $\mathcal{M}^{\perp}$. Furthermore, in \cite{LOYZ17}, an example of a weakly mixing system is presented where the set of distal points is not dense, yet its induced system has dense distal points. By combining the results from \cite{LYY15} and \cite{HSY20}, it is established that a transitive system is disjoint from all minimal systems if and only if its hyperspace system is (\cite[Corollary 4.8]{HSY20}).

In \cite{Opr19}, Oprocha proposed an sufficient condition for a system to be disjoint from all minimal systems. Specifically, he stated that if a system $(X, T)$ is weakly mixing and for every minimal system $(Y, S)$, there exists a countable dense set $D_{Y} \subset X$ satisfying the following condition: for every nonempty open set $U$ in $X$, any point $y \in Y$, and any open neighborhood $V$ of $y$, there exists some $x \in D_{Y} \cap U$ such that the set
$$N_{T \times S}((x,y), U \times V) = \{n \in \mathbb{Z} : T^n x \in U, S^n y \in V\}$$
is syndetic, then $(X,T)$ is disjoint from any minimal system. Later, in \cite{HSY20}, Huang, Shao, and Ye showed that Oprocha's condition is not only sufficient but also necessary for a transitive system to be disjoint from all minimal systems. This finding provides an answer to Furstenberg's problem for transitive systems. However, the characterization mentioned earlier is not entirely intrinsic, as the selection of $D_Y$ is contingent upon the minimal system $(Y,T)$. In this paper, for general topological systems, we will give some more intrinsic characterizations of systems that are disjoint from all minimal systems.

\subsection{Main results}
The following theorem offers a clearer and more intrinsic understanding of the conditions under which a  system is disjoint from all minimal systems.

\begin{thma}[Characterization of systems disjoint from all minimal systems]\label{thm-tran}
Let $(X,T)$ be a topological dynamical system. Then $X \perp \mathcal{M}$  if and only if there exist minimal subsets $\{M_i\}_{i \in \mathbb{N}}$ in $X$  with $\overline{\bigcup_{i \in \mathbb{N}} M_i} = X$  such that each $M_i$ is disjoint from $X$.
\end{thma}

 A measure theoretical analogue of Theorem \ref{thm-tran} is as follows: an invertible probability measure preserving system is disjoint from all ergodic automorphisms if and only if  it is disjoint from almost all its ergodic components.
See Appendix for a proof of this result.

\medskip

For a topological dynamical system $(X,T)$, we say $(X,T)$ having \emph{topological decomposition}, if its maximal equicontinuous factor map $\pi: (X,T) \rightarrow (X_{eq},T_{eq})$ is an open map, and $X_{eq}$ is composed solely of fixed points. (See Subsection \ref{subsec-RP-relation} for the definition of the maximal equicontinuous factor.) The concept of topological decomposition serves as a topological counterpart to the ergodic decomposition in measure-preserving systems.

Inspired by the work of  G\'{o}rska, Lema\'{n}czyk, and de la Rue and the alternative approach by Glasner and Weiss, we present the following characterization for general systems.
As usual the topological analogue of a measure theoretical result is complicated by the necessity to pass to almost one-to-one or even proximal extensions.

\begin{thma}[Characterization of systems disjoint from all minimal systems]\label{thm general}
Let $(X,T)$ be a topological dynamical system. Then  $X\perp \mathcal{M}$ if and only if
$(X,T)$ has an almost one-to-one extension $(X^{*},T^*)$ such that the maximal equicontinous factor map $\pi: X^{*}\rightarrow X^{*}_{eq}$ is a topological decomposition, and there are minimal sets $\{M_i\}_{i\in\mathbb{N}}$ in $(X^{*},T^{*})$ with $\overline{\bigcup_{i=1}^{\infty}M_i}=X^{*}$ and a dense $G_{\delta}$ set $Z^{*}$ of $X_{eq}^{*}\times X_{eq}^{*}$ such that for each $(z_1,z_2)\in Z^{*}$,
\begin{itemize}
\item[(a)]  $(\pi^{-1}(z_1), T^{*})\curlywedge (\pi^{-1}(z_2), T^{*})$, i.e. $(\pi^{-1}(z_1)\times \pi^{-1}(z_2), T^*\times T^*)$ is  transitive;
\item[(b)] $\pi^{-1}(z_i)\perp M_j$ for any $i=1,2$ and $j\in\mathbb{N}$.
\end{itemize}
\end{thma}

In Example \ref{exam-3} of the following text, we will illustrate that the almost one-to-one condition in the theorem is necessary.

A system is {\em semi-simple} if every point within it is a minimal point. We now present the following characterization for semi-simple systems that are disjoint from all minimal systems.

\begin{thma}[Characterization of semi-simple systems disjoint from all minimal systems]\label{thm-ss}
Let $(X,T)$ be a semi-simple system. Then $X\bot\mathcal{M}$ if and only if  the set
\[\Delta_{\perp}(X):= \left\{(x_1,x_2)\in X\times X: (\overline{\O}(x_1,T),T)\perp (\overline{\O}(x_2,T),T)\right\}\]
is residual in $X\times X$.
\end{thma}

For a distal system $(X,T)$, this characterization is also equivalent to certain conditions related to the regionally proximal relation $Q(X)$ of $(X,T)$ 
(for definitions, see subsection \ref{subsec-RP-relation}).

\begin{thma}\label{Main-thm-B}
Let $(X,T)$ be a distal system. Then the following statements are equivalent:
\begin{enumerate}
\item $X\perp\mathcal{M}$.
\item For each minimal subsystem $W$ of $X$ we have $W\times W\subset Q(X)$.
\item There are minimal subsystems $\{X_i\}_{i\in\mathbb{N}}$ with $\bigcup_{i=1}^{\infty} X_i$ is dense in $X$, and for each $i\in\N$, $X_i\times X_i\subset Q(X)$.
\end{enumerate}
\end{thma}

See Section \ref{sec-distal} for more characterizations about distal systems disjoint from all minimal systems.


A key ingredient in the proof of our main theorems is the following property of quasifactors, which may be of independent interest. Any minimal subsystem of $(2^X,T)$ is also called a {\em quasifactor} of $(X,T)$ (see Subsection \ref{subsec-qf} for details).

\begin{thma}[Countability of pairwise disjoint quasifactors]\label{thm QF}
Let $(X,T)$ be a minimal system. Then any collection of pairwise disjoint non-trivial quasifactors of $(X,T)$ is at most countable.
\end{thma}

The countability result in Theorem \ref{thm QF} is inferred from the separability of $C(X)$, the Banach space of all real valued continuous functions on $X$ with supernorm. This highlights the importance of metrizability in this context. The universal minimal flow serves as a counterexample to illustrate this point. Several enhanced versions of Theorem \ref{thm QF} are necessary for proving the general cases, and these are thoroughly discussed in Section \ref{sec QF}. 
It is worth to mention that in the proof given by Glasner and Weiss for ergodic case, the countability of the bases of $L^2(\mu)$ is used.

\medskip

In conclusion, we note that many previous results on disjointness can be extended to systems under group actions. For instance, in \cite{GTWZ21}, the authors demonstrate that the Bernoulli shift of any discrete group $G$ is disjoint from any minimal $G$-system (for a concise proof, see \cite{Ber20}). In \cite{XY}, numerous results originally obtained for $\mathbb{Z}$-actions in \cite{HY05, HSY20} are shown to apply to countable group actions. However, in this paper, we will not explore systems under general group actions. Finally, we remark that our results also hold for continuous maps.

\subsection{Some examples}
Below, we will give some examples to illustrate the theorems mentioned above.

\begin{exam}
Let $X=\T^2=\R^2/\Z^2$ and $T: X\rightarrow X, (x,y)\mapsto (x,x+y) \pmod{\Z^2}$. Then $(X,T)$ is semi-simple. Let $T_x: \{x\}\times \T\rightarrow \{x\}\times \T, (x,y)\mapsto (x,x+y \pmod \Z)$. Then for $x\not\in \Q$, $(\{x\}\times \T,T_x)$ is minimal. And if $x_1,x_2\in \R$ are rationally independent, $(\{x_1\}\times \T,T_{x_1})\perp (\{x_2\}\times \T,T_{x_2})$. Thus by Theorem \ref{thm-ss}, $(X,T)\perp \mathcal M$. Since $(X,T)$ is also distal, we can get $(X,T)\perp \mathcal M$ by Theorem \ref{Main-thm-B} since for each minimal subsystem $W$ of $X$ we have $W\times W\subset Q(X)$.

For a similar 
but more complicated example using Heisenberg group, see \cite[Subsection 3.7]{GLR24}.
\end{exam}

\begin{exam}\label{exam-2}
Furstenberg showed in \cite{Fur67} that a totally transitive system with dense periodic points is disjoint from all minimal systems. Now this result is an easy corollary of Theorem \ref{thm-tran}, since a totally transitive system is disjoint from any periodic system.
In particular, the Bernoulli shift is an element of $\mathcal{M}^{\perp}$.
\end{exam}

In the following example, we will illustrate that the almost one-to-one condition in the theorem \ref{thm general} is necessary.

\begin{exam}\label{exam-3}
Let $(\Sigma_2=\{0,1\}^\Z, \sigma)$ be the full shift. It is a Bernoulli system, and by Example \ref{exam-2}, $(\Sigma_2,\sigma)\perp \mathcal M$. Suppose that the dynamical system $(X^*,T^*)$ is composed of the disjoint union of two Bernoulli systems $\Sigma_2$. Select the same fixed point in each of the two Bernoulli systems, and then identify these two fixed points as a single point to obtain the system $(X,T)$. Then $(X^*,T^*)$ is an almost one-to-one extension of $(X,T)$. It is easy to see that $(X^*,T^*), (X,T)\in {\mathcal M}^\perp$. Since $(\Sigma_2,\sigma)$ is weakly mixing, $X_{eq}^*$ is composed of two fixed point, and $\pi^*: X^*\rightarrow X_{eq}^*$ is a topological decomposition satisfying conditions (a) and (b) in Theorem \ref{thm general}. But for $(X,T)$, $X_{eq}$ is a fixed point. Since $(X,T)$ is not transitive, $\pi: X\rightarrow X_{eq}$ does not satisfy (a) in Theorem \ref{thm general}.

Moreover, there is a distal system disjoint from all minimal system such that $\pi:X\ra X_{eq}$ is not open. Let $X=\T^1\times [-1,1]$ and define $T:X\ra X$ such that $T(x,y)=(x,y)$ for any $(x,y)\in \T^1\times [-1,0)$, and $(x,y)=(x+y,y)$ for $(x,y)\in \T^1\times [0,1]$. This is a distal system and by Theorem \ref{Main-thm-B} it is disjoint from all minimal systems. But, $\pi:X\ra X_{eq}$ is not open. To see this we just note that
$X_{eq}$ can be viewed as the union of a cone and a segment such that the vertex of the cone and one end point of the segment are identified as $p$.
For any $x\in \pi^{-1}p$ and a small open neighborhood $U$, it is clear that $\pi(U)$ is not open in $X_{eq}$.
\end{exam}

\subsection*{Organization of the paper}

The paper is organized as follows: In Section \ref{sec-intro}, we introduce our motivation and present our main results. Section \ref{sec-pre} is dedicated to outlining key concepts and fundamental lemmas that will be utilized throughout the paper.

Section \ref{sec-qf1} presents several foundational lemmas related to disjointness and quasifactors. In Section \ref{sec-residual}, we explore residual properties concerning disjointness. In particular, we show that if $(X,T)$ is a system and $(Y, S)$ is a minimal system, then $(X,T)\perp (Y,S)$ if and only if the set
$\{x\in X: (\overline{\O}(x,T),T)\perp (Y,S)\}$
is residual in $X$.
In Section \ref{sec QF}, our objective is to establish that for any minimal system, there can be at most countably many pairwise disjoint quasifactors. We also present several strengthened versions of this assertion. In Section \ref{p-thmA} we give the proof of Theorem \ref{thm-tran}. In Section \ref{sec-gen}, we provide the proofs for Theorem \ref{thm general} and Theorem \ref{thm-ss}. Section \ref{sec-distal} is dedicated to characterizing distal systems that are disjoint from all minimal systems. Finally, in Section \ref{sec-open}, we pose several open questions.


\section{Preliminaries}\label{sec-pre}

In this section we give some necessary notions and some known facts used in the paper. For more notions and results in topological dynamical systems, one may refer to \cite{Au88, G76}. In the article, integers, nonnegative integers and natural numbers
are denoted by $\Z$, $\Z_+$ and $\N$ respectively.

\subsection{Basics of topological dynamical systems}

\emph{A topological dynamical system} (simply referred to as \emph{a system}) is a pair $(X, T)$, where $X$ is a compact metric space and $T$ is a homeomorphism on $X$. Throughout this paper, we consistently use $\rho_X$ (or simply $\rho$ when there is no risk of ambiguity) to denote the metric on $X$.

Let $(X, T)$ be a system and $x\in X$. Then $\O(x,T)=\{T^nx: n\in \Z\}$ denotes the
{\em orbit} of $x$. We use $\overline{\O}(x,T)$ to denote the closure of the orbit of $x$. 
A subset $A\subseteq X$ is called {\em invariant} (or {$T$-invariant}) if $TA= A$. When $Y\subseteq X$ is a closed and
invariant subset of the system $(X, T)$, we say that the system
$(Y, T|_Y)$ is a {\em subsystem} of $(X, T)$. Usually we will omit the subscript, and denote $(Y, T|_Y)$ by $(Y,T)$.
If $(X, T)$ and $(Y, S)$ are two systems, their {\em product system} is the
system $(X \times Y, T\times S)$ where $T\times S: X\times Y\rightarrow X\times Y, (x,y)\mapsto (Tx,Sy)$.


Let $(X,T)$ be a system. For 
$x,y\in X$, if $\inf_{n\in\mathbb{Z}} \rho(T^{n}x,T^{n}y)>0$ then we say they are \emph{distal}, otherwise we say they are \emph{proximal}.
Denote by ${\bf P}(X,T)$ the set of all proximal pairs of $(X,T)$. If every two distinct points in $X$ are distal, then we say $(X,T)$ is a \emph{distal system}. A system $(X,T)$ is distal if and only if ${\bf P}(X,T)= \Delta(X)=\{(x,x):x\in X\}$.

A system $(X,T)$ is \emph{equicontinuous} if for any $\epsilon>0$ there is some $\delta>0$ such that for any $x,y\in X$, one has $\sup_{n\in\mathbb{Z}}\rho(T^{n}x, T^{n}y)<\epsilon$ whenever $\rho(x,y)<\delta$. Each equicontinuous system is distal.


\subsection{Factor maps}

A {\em factor map} $\pi: X\rightarrow Y$ between two systems $(X,T)$
and $(Y, S)$ is a continuous surjective map which intertwines the
actions (i.e. $\pi\circ T= S\circ \pi$); one says that $(Y, S)$ is a {\it factor} of $(X,T)$ and
that $(X,T)$ is an {\it extension} of $(Y,S)$. The systems are said to be {\em isomorphic} if $\pi$ is bijective.
In this paper, if $(Y,S)$ is a factor of $(X,T)$, it is convenient to use one symbol $T$ to denote the transformation when no room for confusing.
Let $\pi: (X,T)\rightarrow (Y, T)$ be a factor map. Then
$$R_\pi=\{(x_1,x_2):\pi(x_1)=\pi(x_2)\}$$
is a closed invariant equivalence relation, and $Y=X/ R_\pi$.

Let $(X,T)$ and $(Y,T)$ be system and let $\pi: (X,T) \to (Y,T)$ be a factor map.
One says that:
\begin{itemize}
  \item $\pi$ is an {\it open} extension if it sends open sets to open sets; 
  \item $\pi$ is a {\it proximal} extension if
$\pi(x_1)=\pi(x_2)$ implies $(x_1,x_2) \in {\bf P} (X,T)$;
  \item $\pi$ is  a {\it distal} extension if $\pi(x_1)=\pi(x_2)$ and $x_1\neq x_2$ implies $(x_1,x_2) \not\in {\bf P} (X,T)$;
  \item $\pi$ is an {\it almost one-to-one} extension  if there exists a dense $G_\d$ set $X_0\subseteq X$ such that $\pi^{-1}(\{\pi(x)\})=\{x\}$ for any $x\in X_0$;
\item $\pi$ is an {\it equicontinuous} extension if for any $\ep >0$ there exists $\d>0$
such that $\pi(x_1)=\pi(x_2)$ and $\rho(x_1,x_2)<\d$ imply $\rho(T^n x_1,T^n x_2)<\ep$ for any $n\in \Z$.

\end{itemize}

\subsection{Equicontinuous relation and regionally proximal relation}\label{subsec-RP-relation}

Let $(X,T)$ be a system. There is a smallest closed invariant equivalence relation $S_{eq}$ such that the quotient flow $(X/S_{eq},T)$ is equicontinuous \cite[Theorem 1]{EG60}. The equivalence relation $S_{eq}$ is called the {\em equicontinuous structure relation} and the factor $(X_{eq}=X/S_{eq}, T)$ is called the {\em maximal equicontinuous factor} of $(X,T)$.
Note that the maximal equicontinuous factor means if a system $(Z, R)$ is an equicontinuous factor of $(X,T)$ then it is also a factor of $(X_{eq},T)$.

Let $(X,T)$ be a system. For $x,y\in X$, we say there are {\it regionally proximal} if there are subsequences $\{x_i\}_{i\in \mathbb{N}}, \{y_i\}_{i\in \mathbb{N}}$ in $X$ and $\{n_i\}_{i\in \mathbb{N}}$ in $\mathbb{Z}$ such that
\[ x_i\rightarrow x,\, y_i\rightarrow y \ \ \text{and }\ \ \rho(T^{n_i}x_i, T^{n_i}y_i)\rightarrow 0, \text{ as } i\rightarrow \infty.\]
Let $Q(X)=\{(x,y)\in X\times X: x, y \text{ are regionally proximal}\}$. It is easy to see that $\Delta(X)\subseteq Q(X)$ and $Q(X)$ is a symmetric, invariant closed relation on $X$. But, in general, it is not an equivalence relation. It is well-known that (see for example \cite[Chapter 9]{Au88})
\begin{itemize}
\item if let $\mathcal{A}(Q(X))$ be the minimal  closed invariant equivalence relation  of $X$ containing $Q(X)$, then $\mathcal{A}(Q(X))=S_{eq}(X)$ and $X_{eq}=X/\mathcal{A}(Q(X))$. In particular, $(X,T)$ is equicontinuous if and only if $Q(X)=\D(X)$.
\item if $(X,T)$ is minimal then $Q(X)=S_{eq}$ is an equivalence relation and the quotient $X/Q(X)$ is just the maximal equicontinuous factor of $X$.
\end{itemize}

Let $\pi: (X, T)\rightarrow (Y, S)$ be an extension of systems. We say $(y_1,y_2)\in Q(Y)$ has a \emph{$\pi$-lift} if there are $x_i\in \pi^{-1}(y_i), i=1,2$ such that $(x_1,x_2)\in Q(X)$. It is well known that $Q(Y)$ has $\pi$-lifting property if $(X,T)$ is minimal \cite[Theorem 7.8]{Au88}.

The following example illustrates that not every $(x, y) \in Q(Y)$ has a lift, even when dealing with an almost one-to-one extension between two distal systems.

\begin{exam}
There is a counterexample of the almost one-to-one extension $\pi:(X,T)\rightarrow (Y,S)$ between two distal system such that $Q(Y)$ does not $\pi$-lift $Q(X)$.

Let $I=[-1,0]\cup [1,2]$ and $X=I\times \T$. Consider the transformation $T: X \rightarrow X$ defined by $T(x,y)=(x,x+y \pmod{\Z})$.
Consequently, $(X,T)$ forms a distal system. Select points  $u,z_1\in \{0\}\times \T$ and $v,z_2\in \{1\}\times \T$ with $u\neq z_1$ and $v\neq z_2$. By collapsing $z_1,z_2$ into a single point $z$,  we construct another distal system $(Y,S)$. Let $\pi:(X,T)\ra (Y,S)$ be the associated factor map, which constitutes an almost one-to-one extension between the two distal systems.

It is clear that $(\pi(u),\pi(v))\in Q(Y)$, $\pi^{-1}(\pi(u))=\{u\}$, $\pi^{-1}(\pi(v))=\{v\}$ and $(u,v)\not\in Q(X)$.
Therefore, $(\pi(u),\pi(v))$ lacks a $\pi$-lift.
\end{exam}

\subsection{Some facts about hyperspaces}\label{subsec-qf}
Let $X$ be a compact metric space. Let $2^X$ be the collection of nonempty closed subsets of $X$.
Let $\rho$ be the metric on $X$,
then one may define a metric on $2^X$ as follows:
\begin{equation*}
 H_\rho(A,B) = \inf \{\ep>0: A\subseteq B_\ep(B), B\subseteq B_\ep(A)\}
\end{equation*}
where $B_\ep (A)=\{x\in X: \rho(x, A)<\ep\}$.
The metric $H_\rho$ is called the {\em Hausdorff metric} of $2^X$, and $2^X$ is called the {\em hyperspace} of $X$.

\medskip

Let $(X,T)$ be a system. We can induce a system on $2^X$. The
action of $T$ on $2^X$ is given by $TA=\{Ta:a\in A\}$ for each $A \in 2^X$. Then
$(2^X,T)$ is a system and it is called the {\em hyperspace system} or the {\em induced system} of $(X,T)$. Any minimal subsystem of $(2^X,T)$ is also called a {\em quasifactor} of $(X,T)$. The fixed point $\{X\}$ is called the {\em trivial quasifactor}.

\medskip

Let $\{A_i\}_{i=1}^\infty$ be an arbitrary sequence of subsets of $X$. Define
$$\liminf A_i=\{x\in X: \text{for any neighbourhood $U$ of $x$, $U\cap A_i\neq \emptyset$ for all but finitely many $i$}\};$$
$$\limsup A_i=\{x\in X: \text{for any neighbourhood $U$ of $x$, $U\cap A_i\neq \emptyset$ for infinitely many $i$}\}.$$
We say that $\{A_i\}_{i=1}^\infty$ converges to $A$, denoted by $\lim_{i\to \infty} A_i=A$, if
$$\liminf A_i=\limsup A_i=A.$$
Now let $\{A_i\}_{i=1}^\infty\subseteq 2^X$ and $A\in 2^X$. Then $\lim_{i\to\infty} A_i=A$ if and only if $\{A_i\}_{i=1}^\infty $ converges to $A$ in $2^X$ with respect to the Hausdorff metric.

Let $X,Y$ be two compact metric spaces. Let $F: Y\rightarrow 2^X$ be a map and $y\in Y$.
We say that $F$ is {\em upper semi-continuous (u.s.c.)} at $y$ if whenever $\lim y_i=y$, one has that $\limsup F(y_i)\subseteq F(y)$.
We say $F$ is {\em lower semi-continuous (l.s.c.)} at $y$ if whenever $\lim y_i=y$, one has that $\liminf F(y_i)\supset F(y)$. If $F$ is u.s.c. (l.s.c.) at every point of $Y$, then we say that $F$ is u.s.c. (l.s.c.).


It is easy to verify that $F: Y\rightarrow 2^X$ is u.s.c. at $y\in Y$ if and only if for each $\ep>0$ there exists a neighbourhood $U$ of $y$ such that $F(U)\subseteq B_\ep(F(y))$; and $F: Y\rightarrow 2^X$ is l.s.c. at $y\in Y$ if and only if for each $\ep>0$ there exists a neighbourhood $U$ of $y$ such that $F(y)\subseteq B_\ep(F(y'))$ for all $y'\in U$.

Now we give some important examples.

\begin{itemize}
  \item If $f: X\rightarrow Y$ is a continuous surjective map, then it is easy to verify that
$$F=f^{-1}: Y\rightarrow 2^X, y\mapsto f^{-1}(y)$$ is u.s.c.

  \item Let $J\subseteq X\times Y$ be a closed subset that projects onto $ Y$. Then
  $$\varphi_J: Y\rightarrow 2^X, y\mapsto J[y]:=\{x\in X: (x,y)\in J\}$$
  is u.s.c.

  \item Let $(X,T)$ be a system. Then the map $$\overline{\O}: X\rightarrow 2^X, x\mapsto \overline{\O(x,T)}$$ is l.s.c. (see for example \cite[Section 4]{G94}).
\end{itemize}

\medskip

We have the following well known result, for a proof see \cite[p.70-71]{Kura2} and \cite[p.394]{Kura1}, or \cite{Fort}.
\begin{thm}\label{continuous-point-K68}
Let $X,Y$ be compact metric spaces. If $F: Y\rightarrow 2^X$ is u.s.c. (or l.s.c.),
then the points of continuity of $F$ form a dense $G_\delta$ set in $Y$.
\end{thm}


The following fact is easy to prove.
\begin{lem}\label{continuous of orb in equi}
Let $(X,T)$ be a system. If $(X,T)$ is equicontinuous, then the map $\overline{\O} :X\ra 2^X$, $x\mapsto \overline{\O}(x,T)$ is continuous.
\end{lem}
\begin{proof} Since $(X,T)$ is equicontinuous, for each $\ep>0$ there is $\delta>0$ such that if $\rho(x,y)<\delta$ then $\rho(T^nx,T^ny)<\ep$
for any $n\in\Z$. Fix $x\in X$, then for any $y\in B_\delta(x)$, we have $\rho(T^nx,T^ny)<\ep$ for any $n\in\Z$ which implies that
$\rho_H(\overline{\O}(x,T), \overline{\O}(y,T))\le \ep$, since $\overline{\O}(x,T)$ and $\overline{\O}(y,T)$ are minimal.
\end{proof}

\begin{lem}[{\cite[Proposition 7]{BS75}}]\label{QF of equi}
Let $(X,T)$ be a system. If $(X,T)$ is equicontinuous then so is $(2^{X},T)$.
\end{lem}

For a proof of the first statement of the following lemma, see \cite[Theorem 3]{BS75}, and for the second statement see \cite[Theorem 2.5]{G75} or \cite[Corollary 11.19]{Au88}.

\begin{lem}\label{QF of distal}
Let $(X,T)$ be a system. The distality of $(X,T)$ does not necessarily imply the distality of $(2^X,T)$.

If $(X,T)$ is distal and minimal, then any minimal subsystem of $(2^{X},T)$ is distal.
\end{lem}

\subsection{Semi-openness and almost one-to-one extension}
A continuous map $f: X\rightarrow Y$ between topological spaces is \emph{semi-open} if for each nonempty open set $V$ in $X$, the interior of $f(V)$ in $Y$ is not empty.
The following result is well-known, and its proof is straightforward (for example, see \cite[Lemma 2.1]{G07}).
 \begin{lem}\label{residual under semiopen}
 Let $\pi: X\rightarrow Y$ be a continuous surjective map between compact metric spaces. If $\pi$ is semi-open, then for any residual set $S$ in $Y$, $\pi^{-1}(S)$ is residual in $X$.
 \end{lem}

The following result is well-known, see \cite[Proposition 3.1]{V70} or \cite[Theorem 1.2]{Aki04}.
\begin{thm}\label{veech}
Let $\pi: X\rightarrow Y$ be a continuous surjective map between compact metric spaces. Suppose that $\pi$ is semi-open and $R$ is a residual subset of $X$. Then the set
\[ A=\{y\in Y: \pi^{-1}(y)\cap R \text{ is residual in } \pi^{-1}(y)\}\]
is residual in $Y$. Here, being residual in $\pi^{-1}(y)$ refers to the relative topology. In particular, $\pi(R)$ is residual in $Y$.
\end{thm}

The following results are derived directly from the definition of almost one-to-one extension.

\begin{lem} \label{almost1-1}
Let $\pi:(X,T)\rightarrow (Y,S)$ be an almost one-to-one extension between two systems.
\begin{enumerate}
\item If $(Y,S)$ is transitive, then $(X,T)$ is transitive.
\item If $(Y,S)$ is minimal, then $(X,T)$ is minimal and $\pi$ is proximal.
\end{enumerate}
\end{lem}

\begin{lem}[{\cite[Lemma 1.1~(b)]{AG01}}]\label{semi-open for almost 1-1}
Let $\pi: (X, T)\rightarrow (Y, S)$ be an almost one-to-one extension. Then $\pi$ is semi-open.
\end{lem}

\section{Joinings and quasifactors}\label{sec-qf1}

In this section, we discuss joinings and subsystems of hyperspace systems.

\subsection{Joinings and disjointness}

Let $(X, T)$ and $(Y,S)$ be two systems. Recall that a \emph{joining} of the systems $(X, T)$ and $(Y, S)$ is defined as a closed subset of $X \times Y$ that is invariant under the action of $T \times S$ and projects onto both coordinates, $X$ and $Y$.
Denote by $Join(X,Y)$ the set of all joinings of $(X,T)$ and $(Y,S)$. It is clear that $Join(X,Y)$ is a non-empty closed subset of $2^{X\times Y}$
with the Hausforff topology.

We say $(X, T)$ and $(Y,S)$ are \emph{disjoint} if $X\times Y$ is the unique joining of them, i.e. $Join(X,Y)=\{X\times Y\}$, and in this case we denote $(X,T)\perp (Y,S)$ or $X\perp Y$. If $(X\times Y, T\times S)$ is transitive, then we say they are {\em weakly disjoint}, and denote it by $(X,T)\curlywedge (Y,S)$ or $X\curlywedge Y$.

Let $\mathcal{M}$ be the set of all minimal systems. For a system $(X,T)$, we use $X\perp \mathcal{M}$ to denote that $X$ is disjoint from all minimal systems. Let $\mathcal{M}^{\perp}$ be the class of systems that are disjoint from any minimal system.
It is clearly that if each point in  a system $(X,T)$ is a fixed point then  $X\perp \mathcal{M}$.


We list some basic properties of disjointness.

\begin{prop}[\cite{Fur67}]\label{easy}
Let $(X,T)$, $(Y,S)$, $(X_1,T_1)$, $(Y_1,S_1)$ and $(Z,R)$ be systems.
\begin{enumerate}
\item If $(X,T)\perp (Y,S)$, then either $(X, T)$ or $(Y,S)$ is minimal.
\item Let $\pi: X\rightarrow X_1$, $\phi: Y\rightarrow Y_1$ be factor maps. If $X\perp Y$, then $X_1\perp Y_1$. In particular, if $(X,T)\bot (Y,S)$ and $(Z,R)$ is a factor of $(X,T)$, then $(Z,R)\bot (Y,S)$.
\item Two disjoint systems cannot have any common nontrivial factors.
\item If $(X,T)$ and $(Y,S)$ are minimal then  $X\bot Y$ if and only if $(X\times Y,T\times S)$ is minimal.
\end{enumerate}
\end{prop}

The following result is easy to be verified, see for example in \cite[Theorem 2.6-(b)]{AG01}.
\begin{lem}\label{disjoint for almost 1-1}
Let $(X',T')\rightarrow (X,T)$ and $(Y',S')\rightarrow (Y,S)$ be almost one-to-one extensions between two systems.
Then $X'\perp Y'$ if and only if $X\perp Y$.
\end{lem}

\subsection{DDMS-property}

Huang and Ye showed each system in ${\mathcal M}^\perp$ has a dense set of minimal points.

\begin{thm}[{\cite[Theorem 4.3, Theorem 2.6]{HY05}}]\label{dense-mi}
Let $(X,T)$ be a system. Assume that $(X,T)\perp \mathcal{M}$. Then the set of minimal points is dense in $X$.

If in addition $(X,T)$ is transitive, then $(X,T)$ is weakly mixing
\end{thm}

By Theorem \ref{dense-mi}, we know that if a system $(X,T)\bot \mathcal{M}$, then $(X,T)$ has the following the {\em DDMS-property} (DDMS is short for ``dense disjoint minimal subsets"):  there are minimal subsystems  $\{M_i\}_{i\in\mathbb{N}}$ of $(X,T)$ such that each $M_i$ is disjoint from $X$ and $\bigcup_{i=1}^\infty M_i$ is dense in $X$.

\medskip

For a system $(X,T)$, a transitive subsystem $(Y, T)$ is {\em maximal} if it is
maximal among all transitive subsystems by the inclusion. By the Zorn Lemma, each transitive subsystem is contained in a maximal transitive subsystem. In a semi-simple system $(X,T)$, the closure of the orbit of each point is both a minimal subset and a maximal transitive subsystem.

\begin{thm}[{\cite[Theorem 4.5]{HY05}}]\label{hy-general}
Let $(X,T)$ be a system and $(Y,S)$ be a non-trivial minimal system. Then $(X,T)\bot (Y,S)$ if and only if there exist countably many maximal transitive subsystems in $(X,T)$ such that there union is dense in $X$ and each of them is disjoint from $Y$.
\end{thm}

\begin{cor}
Let $(X,T)$ be a dynamical system. Then $(X,T)\bot \mathcal{M}$ if and only if for any
minimal system $(Y,S)$ there exist countably many  maximal transitive subsystems (depending on $Y$) such
that their union is dense in $X$ and each of them is disjoint from $Y$.
\end{cor}


For a distal system, we have the following result.

\begin{thm}[{\cite[Theorem 4.9]{HY05}}]\label{mini-equi}
Let $(X,T)$ be a distal system. Then $X\bot \mathcal{M}$ if and only if $(X,T)$ is disjoint from any minimal equicontinuous system.
\end{thm}

\subsection{Subsystems of $(2^X,T)$ and quasifactors}

The following basic result about quasifactors was given by Glasner.

\begin{thm}[{\cite[Theorem 2.3]{G75}}]\label{thm-quasifactor}
Let $(X,T)$ and $(Y, S)$ be two systems. Let $\varphi: Y\rightarrow 2^X$ be an u.s.c. (or a l.s.c.) map such that $\varphi(Sy)=T\varphi(y), \forall y\in Y$. Let $Y_c$ be the set of continuous points of $\varphi$ and
$$W=\overline{\{\big(y,\varphi(y)\big): y\in Y_c\}}\subset Y\times 2^{X}.$$
Let $\pi_1$ and $\pi_2$ be the projections of $Y\times 2^X$ on $Y$ and $2^X$ respectively.
Then $(W,S\times T)$ is a subsystem of $(Y\times 2^{X}, S\times T)$ and $\pi_1: (W,S\times T)\rightarrow (Y,S)$ is almost one-to-one. And $(\X=\pi_2(W),T)$ is a subsystem of $(2^X,T)$.
$$
\xymatrix{
                & (W, S\times T) \ar[dl]_{\pi_1} \ar[dr]^{\pi_2} \\
 (Y,S) & &     (\X,T)        }
$$
If in addition $(Y,S)$ is minimal, $\X$ is a quasifactor of $(X,T)$.
\end{thm}

\begin{rem}
In the statement of \cite[Theorem 2.3]{G75}, $(X,T)$ and $(Y, S)$ are minimal. But the same proof works without this assumption.
See also \cite[Theorem 4.12]{AG01}.
\end{rem}

The name quasifactor comes from the following construction of Glasner \cite{G75}, which shows that for a minimal system, every factor is, up to an almost one-to-one extension, isomorphic to a quasifactor.

\begin{cor}
Let $\pi: (X,T)\rightarrow (Y,S)$ be a factor map of minimal systems. Then there exists a quasifactor $\X$ of $(X,T)$, which is an almost one-to-one extension of $(Y,S)$.
\end{cor}

\begin{proof}
Let $\varphi=\pi^{-1}: Y\rightarrow 2^X$, and it is u.s.c. Since $\pi: X\rightarrow Y$ is a factor map, we have $\pi^{-1}\circ S=T\circ \pi^{-1}$. Thus we can apply Theorem \ref{thm-quasifactor} to get $W$ and $\X$. One can verify that $\pi_2: W\rightarrow \X$ is one-to-one and hence a homeomorphism. Then $\phi=\pi_1\circ \pi_2^{-1}: \X\rightarrow Y$ is almost one-to-one.
$$
\xymatrix{
                & (W, S\times T) \ar[dl]_{\pi_1} \ar[dr]^{\pi_2} \\
 (Y,S) & &     (\X,T)\ar[ll]_{\phi}  .      }
$$

As $(Y,S)$ is minimal, both $W,\X$ are minimal. In particular, $\X$ is a quasifactor of $(X,T)$.
\end{proof}

Let $(X,T)$ and $(Y, S)$ be two systems. Let $J$ be a joining of $(X,T)$ and $(Y, S)$. Then the map $\varphi_J: Y\rightarrow 2^{X}, y\mapsto J[x]:=\{x\in Y: (x,y)\in J\}$ is u.s.c. and $\varphi_J(Sy)=T\varphi_J(y), \forall y\in Y$. Thus by Theorem \ref{thm-quasifactor}, one can define a subsystem $(\X_J,T)$ of $(2^X,T)$. To be precise, we have the following result, which will be used frequently in the subsequent sections.

\begin{cor}\label{almos 1-1 induced by joining}
Let $(X,T)$ and $(Y, S)$ be two systems. Let $J$ be a joining of $(X,T)$ and $(Y, S)$.  Let $Y_c$ be the set of continuous points of $\varphi_J: Y\rightarrow 2^X$ and
$$W_J=\overline{\{\big(y,\varphi_J(y)\big): y\in Y_c\}}\subset Y\times 2^{X}.$$
Then $(W_J,S\times T)$ is a subsystem of $(Y\times 2^{X}, S\times T)$ and $\pi_1: (W_{J},S\times T)\rightarrow (Y,S)$ is almost one-to-one. And $(\X_J,T)$ is a subsystem of $(2^X,T)$, where $\X_J=\pi_2(W_{J})$.
$$
\xymatrix{
                & (W_J, S\times T) \ar[dl]_{\pi_1} \ar[dr]^{\pi_2} \\
 (Y,S) & &     (\X_J,T)        }
$$

If in addition $(Y,S)$ is minimal, $\X_J$ is a quasifactor of $(X,T)$. In this case, we call $\X_J$ a {\em joining quasifactor}.
$\X_J$ is trivial (i.e. $\X_J=\{X\}$) if and only if $J=X\times Y$.
\end{cor}

\begin{cor}[{\cite[Theorem 2.4]{G75}}]\label{joning quasifactor}
Let $(X,T)$ be a system and $(Y, S)$ be a minimal system. If $(X,T)\not\perp (Y,S)$, then there exists an almost one-to-one extension $(Y^*, H)$ of $(Y,S)$, which has a non-trivial quasifactor $(\X_Y,T)$ of $(X,T)$ as its factor.

We call $(\X_Y,T)$ a {\em joining quasifactor of $(X,T)$ with respect to $(Y,S)$}.
\end{cor}

\begin{proof}
Since $X\not\perp Y$, we take a joining $J$ of $(X,T)$ and $(Y,S)$ such that $J\neq X\times Y$.   Then in Corollary \ref{almos 1-1 induced by joining}, $Y^*=W_J$ and $\X_Y=\X_J$ are what we need.
\end{proof}

\begin{rem}
Since there may be lots of joinings between two systems $(X,T)$ and $(Y,S)$, there may be lots of joining quasifactors of $(X,T)$ with respect to $(Y,S)$.
\end{rem}

\begin{lem}\label{joining gives QF}
Let $(X,T)$ be a system.  Let $(Y,S)$ and $(Z,H)$ be two minimal system such that $X\not\perp Y, X\not\perp Z$. If  $Y\perp Z$, then $\mathcal{X}_{Y}\perp \mathcal{X}_{Z}$, where $\X_Y,\X_Z$ are non-trivial quasifactors of $(X,T)$ defined in Corollary \ref{joning quasifactor}.
\end{lem}

\begin{proof}
By Corollary \ref{joning quasifactor}, there are almost one-to-one extensions $\pi: Y^*\rightarrow Y$ and $\phi: Z^*\rightarrow Z$ such that non-trivial quasifactors  $\X_Y,\X_Z$ are factors of $Y^*$ and $Z^*$ respectively.

Then it follows from Lemma \ref{disjoint for almost 1-1} and $Y\perp Z$ that  $Y^*\perp Z^*$ as $\pi$ and $\phi$ are almost one-to-one. And by Lemma \ref{easy}, we have $\mathcal{X}_{Y}\perp \mathcal{X}_{Z}$ since they are factors of ${Y^*}$ and ${Z^*}$, respectively.
\end{proof}

\section{Residual properties about disjointness}\label{sec-residual}
In this section we give some residual properties about disjointness. In particular, we will show that if $(X,T)$ is a system and $(Y, S)$ is a minimal system, then $(X,T)\perp (Y,S)$ if and only if the set
$\{x\in X: (\overline{\O}(x,T),T)\perp (Y,S)\}$
is residual in $X$.

\medskip

First we have the following easy observation.

\begin{lem}\label{heredity of G_{delta}}
Let $X$ be a complete metric space and $X_0$ be a dense $G_{\delta}$ subset of $X$. If $\Omega$ is a dense $G_{\delta}$ set of $X_0$ under the relative topology, then $\Omega$ is also a dense $G_{\delta}$ subset of $X$.
\end{lem}
\begin{proof}
Since $X_0$ is a dense $G_{\delta}$ subset of $X$, there are dense open sets $U_n$ of $X$ such that  $X_0=\cap_{n=1}^{\infty}U_n$. Similarly, there are dense  open sets $V_n$  of $X$ such that $\Omega=\cap_{n=1}^{\infty} (V_n\cap X_0)$. It follows that
\[ \Omega=\bigcap_{m=1}^{\infty}\bigcap_{n=1}^{\infty} (U_m\cap V_n).\]
Thus $\Omega$ is  a dense $G_{\delta}$ subset of $X$.
\end{proof}

\begin{lem}\label{disjoint of orbits}
Let $(X,T)$ and $(Y, S)$ be two systems, and let $X_0$ and $Y_0$ be the sets of continuous points for the maps $\overline{\O}: X\rightarrow 2^X, x \mapsto \overline{\O}(x,T)$ and $\overline{\O}: Y\rightarrow 2^Y, y \mapsto \overline{\O}(y,S)$, respectively. Then the set
\[\left\{ (x, y)\in X_0\times Y_0: (\overline{\O}(x,T),T)\perp (\overline{\O}(y,S),S)\right\}  \]
is a $G_{\delta}$ subset of $X_0\times Y_0$.
\end{lem}

\begin{proof}
Let $H$ denote the Hausdorff metric on the hyperspace $2^{X\times Y}$.
For $(x,y)\in X\times Y$, let $\mathcal{J}_{x,y}$ be the set of joinings of $(\overline{\O}(x,T),T)$ and $(\overline{\O}(y,S),S)$, i.e. $\mathcal{J}_{x,y}=Join(\overline{\O}(x,T), \overline{\O}(y,S))$.

For each $k\in\mathbb{N}$, define
\[\Omega_{k}=\{ (x,y)\in X_0\times Y_0: {\diam}_{H} (\mathcal{J}_{x,y}) \geq 1/k\}.\]
We claim that $\Omega_{k}$ is closed in $X_0\times Y_0$.  First, observe that if $(x_n,y_n)\rightarrow (x_0,y_0)$ in $X_0\times Y_0$ and $J_{n}\rightarrow J$ in $2^{X\times Y}$ with $J_n\in \mathcal{J}_{x_n,y_n}$, then $J\in \mathcal{J}_{x_0,y_0}$. Indeed, since $\lim \overline{\O}(x_n,T)=\overline{\O}(x_0,T)$,  $\lim \overline{\O}(y_n,S)=\overline{\O}(y_0,S)$, $J_n\subset \overline{\O}(x_n,T)\times \overline{\O}(y_n,S)$ and $J_n$ has full projections onto $\overline{\O}(x_n,T)$ and $\overline{\O}(y_n,S)$, it follows that  $J\subset \overline{\O}(x_0,T)\times \overline{\O}(y_0,S)$ and $J$ has full projections onto $\overline{\O}(x_0,T)$ and $\overline{\O}(y_0,S)$. Clearly, $J$ is $T\times S$ invariant. Thus $J\in\mathcal{J}_{x_0,y_0}$.

Now, suppose that $(x_n,y_n)\rightarrow (x_0,y_0)$ in $X_0\times Y_0$ and $J_{n}\rightarrow J, J'_{n}\rightarrow J'$ in $2^{X\times Y}$ with  $J_n, J'_n\in \mathcal{J}_{x_n,y_n}$ and $H(J_n, J'_n)\geq 1/k$. It remains to show that $H(J, J')\geq 1/k$. This follows trivially from the convergence in the Hausdorff metric.

Recall that in a metric space $Z$, if $A\subset B\subset Z$, then  $A$ is closed in $B$ if and only if $A = B \cap \text{cl}_{Z}(A)$, where $\text{cl}_Z(A)$ is the closure of $A$ in $Z$. Thus, for each $k \in \mathbb{N}$, the closedness of $\Omega_k$ in $X_0 \times Y_0$ implies that $\Omega_k = (X_0 \times Y_0) \cap \overline{\Omega_k}$. Hence,
\begin{align*}
 E:&= (X_0\times Y_0)\setminus (\bigcup_{k=1}^{\infty}\Omega_k)= (X_0\times Y_0)\setminus (\bigcup_{k=1}^{\infty}\overline{\Omega_k})\\
 &= (X_0\times Y_0)\cap \bigcap_{k=1}^\infty  ((X\times Y)\setminus \overline{\Omega_k})
 \end{align*}
is a $G_{\delta}$ subset of $X_0\times Y_0$. It is easy to verify that
\[E=\left\{ (x, y)\in X_0\times Y_0: (\overline{\O}(x,T),T)\perp (\overline{\O}(y,S),S)\right\}.\]
Thus it is a  $G_{\delta}$ subset of $X_0\times Y_0$.
\end{proof}

Now we are ready to give the main result of this section.

\begin{thm}\label{dense of disjoint of orbits}
Let $(X,T)$ be a system and $(Y, S)$ be a minimal system.
\begin{enumerate}
\item[(1)] If $(X,T)\perp (Y,S)$ then the set
\[\{x\in X_0: (\overline{\O}(x,T),T)\perp (Y,S)\}\]
is a dense $G_\delta$ subset of $X$, where $X_0$  is the set of continuous points for $\overline{\O}: X\rightarrow 2^X, x \mapsto \overline{\O}(x,T)$ .
\item[(2)] If the set $\Sigma:=\{x\in X: (\overline{\O}(x,T),T)\perp (Y,S)\}$ is dense in $X$, then $(X,T)\perp (Y,S)$.
\end{enumerate}
\end{thm}

\begin{proof}
(1) Assume that $(X,T)\perp (Y,S)$. Following the proof of Lemma \ref{disjoint of orbits}, for each $k \geq 1$, the set
$\Omega_{k} := \{x \in X_0 : \text{diam}_{H}(\mathcal{J}_{x}) \geq \frac{1}{k}\}$
is closed in $X_0$, where $\mathcal{J}_{x}$ is the set of joinings  of $(\overline{\O}(x,T),T)$ and $(Y,S)$, i.e. $\mathcal{J}_{x}=Join(\overline{\O}(x,T), Y)$. Moreover,
$$\left\{x \in X_0 : \overline{\O}(x,T) \perp Y\right\} = X_0 \setminus \left( \bigcup_{k=1}^{\infty} \Omega_k \right)=X_0\cap \bigcap_{k=1}^\infty  (X\setminus \overline{\Omega_k}).$$
is a $G_\delta$ subset of $X$, as $X_0$ is a $G_\delta$ subset of $X$.

To show that $\{x \in X_0 : \overline{\O}(x,T) \perp Y\}$ is dense in $X$, we employ the method of proof used in  \cite[Theorem 4.5]{HY05}. For the sake of completeness, we provide a proof.

Let $U\subseteq X, V\subseteq Y$ be non-empty open subsets, and let
$$X(U, V )= \{x\in X : N_T(x, U)\cap N_S(y, V ) \neq \emptyset \text{ for any } y\in  Y \},$$
where $N_T(x, U) = \{n\in \mathbb{Z}: T^nx\in U\}$ and $N_S(y, V ) = \{n \in \mathbb{Z} : S^ny\in  V \}$.

\noindent {\bf Claim.}\quad {\em Let $U\subseteq X, V\subseteq Y$ be non-empty open subsets. Then the set
$$\displaystyle \widetilde{X}(U, V):=\big( X \setminus \overline{\bigcup_{n\in \mathbb{Z}} T^n U }\big)\cup  X(U, V)$$
is a $T$-invariant dense open subset of $X$.}

\begin{proof}[Proof of Claim]
It is easy to verify that $X(U,V)$ is invariant and open. Thus $\widetilde{X}(U,V)$ is also invariant and open. It suffices to show that $\widetilde{X}(U,V)$ is dense in $X$.

Now we show that $\overline{\bigcup_{n\in \mathbb{Z}} T^n U}\subseteq \overline{X(U,V)}$, and which implies that $\widetilde{X}(U,V)$ is dense. Let $W\subseteq U$ be a nonempty open subset. We claim that $X(W,V)\cap W\neq \emptyset$. If not, for each $x\in W$, there is some $y_x\in Y$ such that $N_T(x,W)\cap N_S(y_x, V)=\emptyset$. Let
$$J=\left(\Big(X\setminus \bigcup_{n\in \Z}T^nW\Big)\times Y \right)\cup \overline{\bigcup_{n\in \Z}\bigcup_{x\in W}(T\times S)^n(x,y_x)}.$$
Then $J$ is a joining of $X$ and $Y$. As $N_T(x,W)\cap N_S(y_x, V)=\emptyset$, for each $x\in W$,
$$\overline{\bigcup_{n\in \Z}\bigcup_{x\in W}(T\times S)^n(x,y_x)}\cap (W\times V)=\emptyset.$$
Thus, $J\cap (W\times V)=\emptyset$. Hence $J\neq X\times Y$, which contradicts $X\perp Y$.
Hence we have that $X(W,V)\cap W\neq \emptyset$ for any non-empty open subset $W$ of $U$. In particular, $X(U,V)\cap W\neq \emptyset$
for any non-empty open subset $W$ of $U$. Since $X(U,V)$ is invariant, we have that $\overline{\bigcup_{n\in \mathbb{Z}} T^n U}\subseteq \overline{X(U,V)}$. The proof of Claim is complete.
\end{proof}

Let $\{V_j\}_{j=1}^{+\infty}$ and $\{U_m\}_{m=1}^\infty$ be basis of $Y$ and $X$ respectively. Put
$$ R_X=\bigcap_{m=1}^{+\infty}\bigcap_{j=1}^{+\infty}\widetilde{X}(U_m, V_j ).$$
By the Claim, $\widetilde{X}(U_m, V_j )$ is a dense open invariant subset of $X$ for $m, j\in \mathbb{N}$, and hence
$R_X$ is a dense $G_\delta$ set of $X$. Now for each $x\in R_X\cap X_0$,  we show that $(\overline{\O}(x,T), T )\perp  (Y,S)$.
Once this is shown, then $\{x \in X_0 : \overline{\O}(x,T) \perp Y\}$ is dense in $X$ as it contains the dense $G_\delta$ set $R_X\cap X_0$.

Let $x\in R_X\cap X_0$. Let $J$ be a joining of $(\overline{\O}(x,T),T)$ and $(Y,S)$.   Take $y\in Y$ such that $(x, y)\in J$.
For any open neighborhood $U$ of $x$ and any non-empty open set $V$ of $Y$, there exist $U_m, V_j$
such that $x\in  U_m\subset U$ and $V_j\subset V$. By the definition of $R_X$, $x\in X(U_m, V_j )$. In particular,
$N_T(x, U_m)\cap  N_S(y, V_j )\neq \emptyset$. As $T\times S(J) = J$ and $(x,y)\in J$, one has $J\cap U_m\times V_j\neq \emptyset$.
Thus $J\cap (U\times V)\neq \emptyset$ for  any open neighborhood $U$ of $x$ and any non-empty open set $V$ of $Y$. It follows that
$J\supset \{x\}\times Y$ and hence $J = \overline{\O}(x,T) \times Y$. Therefore, $(\overline{\O}(x,T), T )\perp  (Y,S)$.

\medskip
(2) Let $J$ be a joining of $(X,T)$ and $(Y,S)$.   Fix an $x\in \Sigma$. Then  $(\overline{\O}(x,T),T)\perp (Y,S)$. Since $J$ is a joining, there is some $y\in Y$ with $(x,y)\in J$. Consequently, $$\overline{\O}((x,y),T\times S)\subset J$$
since $J$ is closed and $T\times S(J)=J$. As $Y$ is minimal, it is clear that $\overline{\O}((x,y),T\times S)$ is a joining of $(\overline{\O}(x,T),T)$ and $(Y,S)$. Then it follows from the disjointness of  $\overline{\O}(x,T)$ and $Y$ that $\overline{\O}((x,y),T\times S)=\overline{\O}(x,T)\times Y$. In particular, $J\supset \{x\}\times Y$.
By the arbitrariness of $x\in \Sigma$, one has $J\supset \Sigma\times Y$.   Since $\Sigma$ is dense in $X$ and $J$ is closed in $X\times Y$, this implies that $J=X\times Y$, and hence $X$ is disjoint from $Y$.
\end{proof}

By Theorem \ref{dense of disjoint of orbits}, one can easily obtain the following.
\begin{cor}\label{many pt for disjoint}
Let $(X, T)$ be a system and $(Y,T)$ be a minimal system. Then $(X,T)\perp (Y,S)$ if and only if  the set
\[\{x\in X: (\overline{\O}(x,T),T)\perp (Y,S)\}\]
is residual in $X$.
\end{cor}

Now we can show one direction of Theorem \ref{thm-ss}.

\begin{cor}\label{dense G of disjoint of orbits}
Let $(X,T)$ be a semi-simple system. If $(X,T)\perp\mathcal{M}$, then the set
\[\Delta_{\perp}(X):= \left\{(x_1,x_2)\in X\times X: (\overline{\O}(x_1,T),T) \perp (\overline{\O}(x_2,T),T)\right\}\]
is residual in $X\times X$.
\end{cor}

\begin{proof}
Let $X_0$  be the set of continuous points for $\overline{\O}: X\rightarrow 2^X, x\mapsto \overline{\O}(x,T)$, which is a dense $G_{\delta}$ subset of $X$.  By Lemma \ref{disjoint of orbits},
$$\Delta_{\perp}(X_0)= \left\{(x_1,x_2)\in X_0\times X_0: (\overline{\O}(x_1,T),T)\perp (\overline{\O}(x_2,T),T)\right\}$$
is a $G_{\delta}$ subset of $X_0\times X_0$. As $X_0$ is a $G_\delta$ subset of $X$, it follows from Lemma \ref{heredity of G_{delta}} that  $\Delta_{\perp}(X_0)$ is a  $G_{\delta}$ subset of $X\times X$.

Next fix
$x_2\in X_0$. As $(X,T)$ is semi-simple, $(\overline{\O}(x_2,T),T)$ is minimal and $X \perp \overline{\O}(x_2,T)$. Thus by Theorem \ref{dense of disjoint of orbits}-(1),
 \[\{x_1\in X_0: (\overline{\O}(x_1,T),T) \perp (\overline{\O}(x_2,T),T) \}\]
is a dense $G_\delta$ subset of $X$. Moreover as $X_0$ is dense in $X$, one has $\Delta_{\perp}(X_0)$ is a dense $G_\delta$ subset of $X\times X$.
Therefore, $\Delta_{\perp}(X)$ is residual in $X\times X$ because $\Delta_{\perp}(X)\supset \Delta_{\perp}(X_0)$.
\end{proof}

We remark here the semi-simplicity in Corollary \ref{dense G of disjoint of orbits} cannot be omitted. For example, consider a transitive non-minimal system in $\mathcal{M}^{\perp}$ such as the Bernoulli shift.

\section{Orders of transitive hyperspace systems and disjoint quasifactors}\label{sec QF}

In this section, our objective is to show that for any minimal system, there can be at most countably many pairwise disjoint non-trivial quasifactors. Additionally, we will present several strengthened versions of this statement. These findings are essential for establishing the main results of the paper.

\subsection{Order of a transitive subsystem of $(2^X,T)$}
We start with an easy fact.

\begin{lem}\label{lem-1}
Let $(X,T)$ be a minimal system. For any subsystem $\mathcal{X}$ of $(2^X,T)$, we have that $\bigcup_{A\in\mathcal{X}} A=X$.
\end{lem}

\begin{proof} 
It suffices to show that for any $x\in X$, we have $x\in \bigcup_{A\in\mathcal{X}}A$. Let $B\in \X$ and $b\in B$. Since $(X,T)$ is minimal, there is some sequence $\{n_i\}_{i\in \N}\subseteq \Z$ such that $T^{n_i}b\to x, i\to\infty$. Without loss of generality, we may assume that $T^{n_i}B\to B', i\to\infty$ in $2^X$ for some $B'\in 2^X$. As $\X$ is invariant, $B'\in \X$. By the definition of $B'$, $x\in B'$ and $x\in \bigcup_{A\in\mathcal{X}}A$.
\end{proof}


\begin{lem}\label{lem-2}
Let $(X,T)$ be a minimal system and $\mathcal{X}, \mathcal{Y}$ be two subsystems of $(2^X,T)$. If $\mathcal{X}\perp \mathcal{Y}$, then for any $A\in \mathcal{X}$ and $B\in \mathcal{Y}$, we have $A\cap B\neq \emptyset$.
\end{lem}

\begin{proof}
Define
\[J=\{(A,B)\in \mathcal{X}\times \mathcal{Y}: A\cap B\neq\emptyset\}.\]
Clearly, $J$ is a closed subset of $\mathcal{X}\times \mathcal{Y}$  and invariant under $T\times T$. Since $(X,T)$ is minimal, by Lemma \ref{lem-1} one has
\[ \bigcup_{A\in\mathcal{X}} A=X, \ \  \bigcup_{B\in\mathcal{Y}} B=X.\]
This implies that for any $A\in \mathcal{X}$ there is some $B\in \mathcal{Y}$ with $A\cap B\neq\emptyset$, and for any $B\in \Y$ there is some $A\in \X$ with $B\cap A\neq \emptyset$. Thus $J$ projects onto $\mathcal{X}$ and  $\mathcal{Y}$ through its  coordinate projections.
Therefore, $J$ is a joining of $\mathcal{X}$ and $\mathcal{Y}$, and hence $J=\mathcal{X}\times\mathcal{Y}$, since  $\mathcal{X}\perp \mathcal{Y}$.  Consequently,  for any $A\in \mathcal{X}$ and $B\in \mathcal{Y}$, we have $A\cap B\neq \emptyset$.
\end{proof}

Now we define the order for a non-empty closed subset in a system.
\begin{de}
Let $(X,T)$ be a system. For $A\in 2^X$, define the {\em order of $A$} as follows:
$$\text{ord}(A)=\sup\ \{n\in \N: A\cap TA\cap \cdots \cap T^{n-1}A \not = \emptyset\}\in \N\cup\{\infty\}.$$
\end{de}

Clearly, $\text{\rm ord}(X)=\infty$. And for a system $(X,T)$ and $x\in X$, if $x$ is a fixed point, ${\rm ord }(\{x\})=\infty$; else ${\rm ord }(\{x\})=1$.

If ${\rm ord}(A)<\infty$, then $\bigcap_{i=0}^{{\rm ord}(A)-1}T^iA \neq \emptyset$ and $A\cap TA\cap \cdots \cap T^{{\rm ord}(A)}A=\emptyset$.

\medskip

Now we give some properties of the map {\rm ord} and define the order for transitive subsystems of hyperspace systems.

\begin{lem}\label{property of order}
Let $(X,T)$ be a system.
\begin{enumerate}
\item The map $\text{\rm ord}: 2^X\rightarrow \N \cup\{\infty\}$ is  upper semi-continuous (u.s.c.).

\item When $(X,T)$ is minimal and $A\in 2^X$ with $A\not=X$, we have ${\rm ord}(A)<\infty$.

\item  If $(X,T)$ is minimal and $\mathcal{W}$ is a non-trivial transitive subsystem of $(2^X,T)$, then there is some $k\in\mathbb{N}$ such that $\text{\rm ord}(A)=k$ for each transitive point $A$ of $\mathcal{W}$ and $\text{\rm ord}(B)\geq k$ for any $B\in \mathcal{W}$.
    In this case, we call $k$ the {\em order of $\mathcal W$}, and denote it by ${\rm ord}({\mathcal W})$.\newline
    In particular, if $\mathcal W$ is a quasifactor, then $\text{\rm ord}(A)=k$ for each point $A$ of $\mathcal{W}$.
\end{enumerate}
\end{lem}

\begin{proof} (1) Let $A_i\in 2^X, i\in \mathbb{N}$ and $A\in 2^X$ satisfy $A_i\ra A$ under the Hausdorff metric. For any given integer $t$ with $0\le t<\limsup_{i\to \infty} {\rm ord}(A_i)$, we can find a subsequence $\{i_1<i_2<\cdots\}$ of nature numbers
such that ${\rm ord}(A_{i_n})>t$ for each $n\in \mathbb{N}$. Thus $\bigcap_{j=0}^t T^j A_{i_n}\neq \emptyset$ for each $n\in \mathbb{N}$. Note that $A_{i_n}\ra A$ when $n\ra \infty$
we have $\bigcap_{j=0}^t T^j A\neq \emptyset$. Thus ${\rm ord}(A)>t$. By the arbitrariness of $t\in [0,\limsup_{i\to \infty} {\rm ord}(A_i))$, we have
${\rm ord}(A)\ge \limsup_{i\to \infty} {\rm ord}(A_i)$. Hence the map $\text{{\rm ord}}: 2^X\rightarrow \N \cup\{\infty\}$ is u.s.c.
\medskip

(2) Suppose  $(X,T)$ is minimal. Note that if $(X,T)$ is a minimal system under $\Z$-action, then both $(X,T)$ and $(X, T^{-1})$ are minimal for the $\Z_+$-actions (see for example \cite[Page 129]{Walters}). Then $\{T^nx:n\in\N\}$ and $\{T^{-n}x:n\in\N\}$ are dense in $X$ for
each $x\in X$. If $\text{{\rm ord}}(A)=\infty$, then $\bigcap_{j=0}^{+\infty}T^{j}A\neq \emptyset$. We take a point $x_0\in \bigcap_{j=0}^{+\infty}T^{j}A$, then $T^{-j}x_0\in A$ for each $j\in\mathbb{N}$. It implies that $A=X$ since $A$ is closed and $\{T^{-n}x_0:n\in\mathbb{N}\}$ is dense in $X$. Hence, when $A\neq X$,  $\text{{\rm ord}}(A)<\infty$.
\medskip

\medskip
(3) Suppose that $(X,T)$ is minimal and $\mathcal{W}$ is a non-trivial transitive subsystem of $(2^X,T)$. Then by (2), $\text{{\rm ord}}(C)<\infty$ for any $C\in \mathcal{W}$.
Now if $A$ is a transitive point of $\mathcal{W}$, then by (1) ${\rm ord}(A)\le {\rm ord}(B)$ for any $B\in \mathcal{W}$ since there is a sequence $\{n_i\}_{i=1}^\infty\subset \Z$
with $T^{n_i}A\to B, i\to\infty$ and ${\rm ord}(T^tA)={\rm ord}(A)$ for any $t\in\mathbb{Z}$. If $A_1,A_2$ are transitive points of $\mathcal{W}$, then it is clear that ${\rm ord}(A_1)={\rm ord}(A_2)\in \mathbb{N}$. Hence there is $k\in\mathbb{N}$ such that $\text{{\rm ord}}(A)=k$ for each transitive point $A$ of $\mathcal{W}$ and $\text{{\rm ord}}(B)\geq k$ for any $B\in \mathcal{W}$.
\end{proof}


\subsection{Minimal subsystems in the hyperspace}

The following result is a generalization of Lemma \ref{lem-2}.

\begin{lem}\label{small disjoint}
Let $(X,T)$ be a minimal system. Suppose that  $\mathcal{X}$ is a subsystem of $(2^{X},T)$ and $\mathcal{Y}$ is a minimal subsystem of $(2^{X},T)$ (i.e. a quasifactor).  If $\mathcal{X}\perp \mathcal{Y}$, then for any $A\in \mathcal{X}$ and $B\in \mathcal{Y}$, we have
\[ A\cap  TB\cap \cdots\cap T^{{\rm ord}(\mathcal{Y})}B \neq\emptyset.\]
\end{lem}

\begin{proof} Without loss of generality, we may assume that the minimal system $\mathcal{Y}$ is non-trivial. Hence ${\rm ord}(\mathcal{Y})=k\in\mathbb{N}$ by Lemma \ref{property of order}~(3).

Consider
\[ J=\{ (A,B)\in \mathcal{X}\times \mathcal{Y}: A\cap  TB\cap \cdots\cap T^{k}B \neq\emptyset\}.\]
We claim that $J$ is a joining of $\mathcal{X}$ and $\mathcal{Y}$.

\medskip
(1)  $J$ is closed in $2^{X}\times 2^{X}$.

Let $(A_n, B_{n})\in J$ with $(A_n,B_n)\rightarrow (A,B)$. Now we show that $(A,B)\in J$. For each $n\in\mathbb{N}$, there exists
\begin{equation}\label{eq1}
 x_n\in A_{n}\cap  TB_{n}\cap \cdots\cap T^{k}B_{n}.
 \end{equation}
We may assume that $x_n\rightarrow x$ in $X$ as $n\to\infty$. By (\ref{eq1}), we have
\[ T^{-1}x_n, \ldots, T^{-k}x_n\in B_{n}.\]
Since $x_n\ra x$, we have $T^{-i}x_n\ra T^{-i}x$ for each $i=1,\ldots,k$. Thus we have
\[ x\in A \text{ and } T^{-i}x\in B, \forall i=1,\ldots, k.\]
In particular,
\[x\in A\cap  TB\cap \cdots\cap T^{k}B\neq\emptyset.\]
Hence $(A,B)\in J$ and $J$ is closed.

\medskip
(2) $J$ is $T\times T$-invariant.

Let $(A,B)\in J$. Then $A\cap  TB\cap \cdots\cap T^{k}B \neq\emptyset$. Since
\[TA\cap T(TB)\cap \cdots\cap T^{k}(TB)=T(A\cap TB\cap \cdots\cap T^{k}B)  \neq\emptyset \]
and
\[T^{-1}A\cap T(T^{-1}B)\cap \cdots\cap T^{k}(T^{-1}B)=T^{-1}(A\cap TB\cap \cdots\cap T^{k}B)  \neq\emptyset, \]
we have $(TA,TB)\in J$ and $(T^{-1}A,T^{-1}B)\in J$. Thus $J$ is invariant.

\medskip
(3) $J$  projects onto each coordinate through its coordinate projections.

First, it follows from the definition of the order that
$$TB\cap \cdots T^k B=T(B\cap \cdots T^{k-1} B)\neq\emptyset $$ for each $B\in\mathcal{Y}$. Meanwhile, by Lemma \ref{lem-1}, $\bigcup_{A\in \mathcal{X}}A=X$. Thus for each $B\in \mathcal{Y}$, there is some $A\in \mathcal{X}$ such that $A\cap TB\cap \cdots\cap T^{k}B\neq\emptyset$. This shows that $J$ projects onto $\mathcal{Y}$ through the second coordinate projection.

Similarly, to show that $J$ projects onto $\mathcal{X}$ through the first coordinate projection, it suffices to show that
\[\bigcup_{B\in\mathcal{Y}}(TB\cap \cdots\cap T^{k}B)=X.\]
By the proof of (1) and (2), $\{ TB\cap \cdots\cap T^{k}B: B\in \Y\}$ is also invariant closed subset of $2^X$. By Lemma \ref{lem-1}, we have $\bigcup_{B\in\mathcal{Y}}(TB\cap \cdots\cap T^{k}B)=X$.


\medskip
Therefore, $J$ is a joining of $\mathcal{X}$ and $\mathcal{Y}$. Hence $J=\mathcal{X}\times \mathcal{Y}$ as $\mathcal{X}\perp \mathcal{Y}$. Then for any $A\in \mathcal{X}$ and $B\in \mathcal{Y}$, we have
\[ A\cap  TB\cap \cdots\cap T^{{\rm ord}(\mathcal{Y})}B \neq\emptyset.\] by the definition of $J$. The proof is complete.
\end{proof}

Now we are ready to prove Theorem \ref{thm QF}.

\begin{proof}[Proof of Theorem \ref{thm QF}]
Let $\Gamma$ be a collection of pairwise disjoint quasifactors of $X$. Without loss of generality, one may assume each quasifactor in $\Gamma$ is nontrivial. Now we show that $\Gamma$ is at most countable.
For each $\mathcal{X}\in \Gamma$, we fix some $A_{\mathcal{X}}\in \mathcal{X}$. Then by Lemma \ref{property of order}~(3), $${\rm ord}(A_{\mathcal{X}})={\rm ord}(\mathcal{X})<\infty.$$ Moreover, it follows from the definition of the order that
\begin{align*}
\bigcap_{j=0}^{{\rm ord}(\mathcal{X})-1}T^j A_{\mathcal{X}}\neq \emptyset \text{ and }A_{\mathcal{X}}\cap TA_{\mathcal{X}}\cap\cdots\cap T^{{\rm ord}(\mathcal{X})}A_{\mathcal{X}}=\emptyset.
\end{align*}
By Urysohn's Lemma, there is some $f_{\mathcal{X}}\in C(X)$ such that
\[f_{\mathcal{X}}\mid_{A_{\mathcal{X}}}\equiv 0\ \ \text{and}\ \ f_{\mathcal{X}}\mid_{TA_{\mathcal{X}}\cap\cdots\cap T^{{\rm ord}(\mathcal{X})}A_{\mathcal{X}}}\equiv 1.\]
However, by Lemma \ref{small disjoint}, for any distinct $\mathcal{X},\mathcal{Y}\in \Gamma$,
\[ A_{\mathcal{X}}\cap TA_{\mathcal{Y}}\cap  \cdots\cap T^{{\rm ord}(\mathcal{Y})}A_{\mathcal{Y}}\neq\emptyset,\]
since $\mathcal{X}\perp\mathcal{Y}$. Hence
\[\|f_{\mathcal{X}}-f_{\mathcal{Y}}\|_{\sup}\geq 1,\]
since 
$f_{\mathcal{X}}(x)=0$ and $f_{\mathcal{Y}}(x)=1$ for any $x\in  A_{\mathcal{X}}\cap TA_{\mathcal{Y}}\cap \cdots\cap T^{{\rm ord}(\mathcal{Y})}A_{\mathcal{Y}}$.

Now, it follows from the separability of $C(X)$ that $\Gamma$ is at most countable.
\end{proof}

\begin{rem}
We note that both metrizability and minimality are essential conditions in Theorem \ref{thm QF} 
The universal minimal system and the Bernoulli shift serve as counterexamples for these conditions, respectively.
\end{rem}

An immediate consequence of Theorem \ref{thm QF} 
is the following interesting result on the disjointness.

\begin{cor}\label{cor-1}
Let $\Gamma=\{(X_{\alpha}, T)\}_{\alpha\in\Lambda}$ be a collection of uncountable pairwise disjoint minimal systems. Then for any minimal system $(Y,S)$ there is some system in $\Gamma$ which is disjoint from $Y$.
\end{cor}

\begin{proof}
Let $(Y,S)$ be a minimal system. Assume $Y\not\perp X_{\alpha}$ for each $\alpha\in\Lambda$. Then, by Lemma \ref{joining gives QF}, this assumption leads to a collection of pairwise disjoint non-trivial quasifactors $\{\mathcal{Y}_{X_\alpha}\}_{\alpha\in\Lambda}$ of $(Y,S)$.
Since the index set $\Lambda$ is an uncountable set, this situation contradicts Theorem \ref{thm QF}, 
which states that any collection of pairwise disjoint quasifactors of a minimal system is at most countable. Thus $Y\perp X_{\alpha}$ for some $\alpha\in\Lambda$. The proof is complete.
\end{proof}

\subsection{Transitive subsystems in the hyperspace}

Now we generalize Corollary \ref{cor-1} to a family of transitive systems. First we have the following lemma which is a generalization of Lemma \ref{small disjoint}.

\begin{lem}\label{lem-tran-1}
Let $(X,T)$ be a minimal system. Suppose that  $\mathcal{X},\Y$ are transitive subsystems of $(2^{X},T)$.
Assume that there exists a family of quasifactors $\{\Y_i\}_{i\in \N}$ such that $\Y\subseteq \overline{\bigcup_{i=1}^\infty \Y_i}$ in $2^X$, ${\rm ord}(\mathcal{Y}_i)\geq {\rm ord}(\mathcal{Y})$ and $\X\perp \Y_i$ for each $i\in\mathbb{N}$.
Then for any $A\in \mathcal{X}$ and $B\in \mathcal{Y}$, we have
\[ A\cap  TB\cap \cdots\cap T^{{\rm ord}(\mathcal{Y})}B \neq\emptyset.\]
\end{lem}

\begin{proof}
Since $\X\perp \Y_i$ for each $i$, by Lemma \ref{small disjoint} we have that for each $A\in \X$ and $C\in \Y_i$,
\[ A\cap TC\cap\cdots\cap T^{{\rm ord}(\mathcal{Y}_{i})}C\neq\emptyset.\]
In particular, since ${\rm ord}(\mathcal{Y}_{i})\geq {\rm ord}(\mathcal{Y})$, it follows that
\begin{equation}\label{eq-1}
A\cap TC\cap\cdots\cap T^{{\rm ord}(\Y)}C\neq\emptyset,
\end{equation}
for any $A\in \X$ and $C\in \bigcup_{i=1}^\infty \mathcal{Y}_i$.

Now since $\Y\subseteq \overline{\bigcup_{i=1}^\infty \Y_i}$, for each $B\in \Y$ there is a sequence $C_{n_i}\in \Y_{n_i}$ with $\{n_i\}_{i\in \N}\subseteq \N$ such that $\lim_{i\to\infty}C_{n_i}=B$ in $2^X$. Meanwhile, each $C_{n_i}$ satisfies (\ref{eq-1}). Taking limit as $i\to \infty$, we conclude that $A \cap TB \cap \cdots \cap T^{{\rm ord}(\Y)}B \neq \emptyset$ for $A\in \X$ and $B\in \Y$. The proof is complete.
\end{proof}

The following result is a generalization of Theorem \ref{thm QF}. 
\begin{thm}\label{lem-tran-2}
Let $(X,T)$ be a minimal system. Let $\{\X_\a\}_{\a\in \Gamma}$ be a collection of transitive subsystems of $(2^X, T)$. Assume that for each $\a\neq \b\in \Gamma$, there exist a family of non-trivial quasifactors $\{\Y_i\}_{i\in \N}$ such that $\X_\b\subseteq \overline{\bigcup_{i=1}^\infty \Y_i}$ in $2^X$, ${\rm ord}(\mathcal{Y}_i)\geq {\rm ord}(\mathcal{X}_\b)$ and $\X_\a\perp \Y_i$ for each $i\in\mathbb{N}$.
Then $\Gamma$ is at most countable.
\end{thm}

\begin{proof}
The proof is almost the same to the one of Theorem \ref{thm QF}. 
For $\a\in \Gamma$, we take a transitive point $A_{\alpha}$ of $(\mathcal{X}_{\alpha},T)$. Then by Lemma \ref{property of order}~(3) and the definition of ${\rm ord}(A_{\alpha})$, we have
$$\bigcap_{j=0}^{{\rm ord}(A_{\alpha})-1}T^jA_{\alpha}\neq \emptyset\text{ and } A_{\alpha}\cap TA_{\alpha}\cap \cdots \cap T^{{\rm ord}(A_{\alpha})}A_{\alpha}=\emptyset.$$ By Urysohn's Lemma, there is some $f_{\alpha}\in C(Y)$ such that
\[ f_{\alpha}\mid_{A_{\alpha}}\equiv 0\ \ \text{and}\ \ f_{\alpha}\mid_{TA_{\alpha}\cap \cdots \cap T^{{\rm ord}(A_{\alpha})}A_{\alpha}}\equiv 1.\]
Now, for any  $\alpha\neq\beta\in \Gamma$, by Lemma \ref{lem-tran-1}, we have
\[ \|f_{\alpha}-f_{\beta}\|_{\sup}\geq 1,\]
because $f_{\alpha}(x)=0$ and $f_{\beta}(x)=1$ for any $x\in  A_{\alpha}\cap TA_{\beta}\cap \cdots\cap T^{{\rm ord}(A_{\b})}A_{\beta}$.
Thus, the separability of $C(X)$ implies the countability of $\Gamma$.
The proof is complete.
\end{proof}

Now we have the following generalization of Corollary \ref{cor-1}.

\begin{prop}\label{disjoint-tran-QF}
Let $\{(X_{\alpha},T_{\alpha})\}_{\alpha\in \Gamma }$ be a collection of transitive systems such that for each  $\alpha\neq \beta\in \Gamma$, there are minimal sets $\{M_i\}_{i\in\mathbb{N}}$ in $X_{\beta}$ with $\overline{\bigcup_{i=1}^{\infty}M_i}=X_{\beta}$ and $X_{\alpha}\perp M_i$ for each $i\in\mathbb{N}$. If $\Gamma$ is uncountable, then for any minimal system $(Y,S)$ there is some $\alpha\in \Gamma$ such that $X_{\alpha}\perp Y$.
\end{prop}

\begin{proof}
Suppose that we can find a minimal system $(Y,S)$ that is not disjoint from each $X_{\alpha}$, $\alpha\in \Gamma$.  We show that $\Gamma$ is countable. To this aim, we will construct a family of transitive subsystems $\{\Y_\a\}_{\a\in \Gamma}$ of $(2^Y,S)$ satisfying conditions in Theorem \ref{lem-tran-2}.

For each $\alpha\in \Gamma$, since $X_{\alpha} \not\perp Y$, by Corollary \ref{almos 1-1 induced by joining}, there is some non-trivial subsystem $(\Y_\a, S)$ of $(2^Y,S)$. To be precise, first choose a joining $J_{\alpha}\in Join(Y,X_\a)$ such that $J_{\alpha}\neq Y\times  X_{\alpha}$. Let $\Omega_{\alpha}$ be the set of continuous points of the mapping $\varphi_\a: X_{\alpha}\rightarrow 2^{Y}, x\mapsto  J_{\alpha}[x]$ and $W_{\alpha}=\overline{\{(x, J_{\alpha}[x])\in X_{\alpha}\times 2^{Y}: x\in\Omega_{\alpha} \}}$. Then the first coordinate projection $\pi^{\alpha}_1: (W_{\alpha},T_{\alpha}\times S)\rightarrow (X_{\alpha},T_{\alpha})$ is an almost one-to-one extension and $\mathcal{Y}_{\alpha}=\pi^{\alpha}_2(W_{\alpha})$ is a nontrivial subsystem $(\Y_\a, S)$ of $(2^Y,S)$. Since $(X_{\alpha},T_{\alpha})$ is transitive, both $W_\a$ and $\mathcal{Y}_{\alpha}$ are transitive.
$$
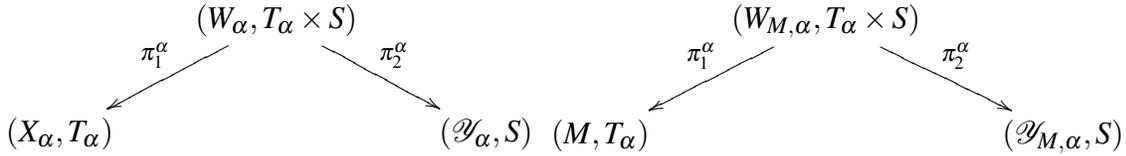
\xymatrix{
                & (W_\a, T_\a\times S) \ar[dl]_{\pi^\a_1} \ar[dr]^{\pi^\a_2} \\
 (X_\a,T_\a) & &     (\Y_\a, S)        }
 \xymatrix{
                & (W_{M,\a}, T_\a\times S) \ar[dl]_{\pi^\a_1} \ar[dr]^{\pi^\a_2} \\
 (M,T_\a) & &     (\Y_{M,\a}, S)        }
$$
Now for each minimal set $M$ in $(X_{\alpha},T_{\alpha})$, we let $\Omega_{M, \alpha}$ be the set of continuous points of the mapping $\varphi_\a|_M: M\rightarrow 2^{Y}, x\mapsto  J_{\alpha}[x]$ and $W_{M, \alpha}=\overline{\{(x, J_{\alpha}[x])\in M\times 2^{Y}: x\in\Omega_{M, \alpha} \}}.$
Then $(W_{M,\alpha}, T_{\alpha}\times S)$ is a minimal system, and $\mathcal{Y}_{M, \alpha}=\pi_2^\alpha(W_{M, \alpha})$ is a minimal subsystem of $(2^Y,S)$, i.e. a quasifactor.

Let $k_{\alpha}={\rm ord}(\mathcal{Y}_{\alpha})$ and take a transitive point $A_{\alpha}=J_{\alpha}[x_{\alpha}]$ of $(\mathcal{Y}_{\alpha},S)$ with $x_{\alpha}\in \Omega_{\alpha}$. Then ${\rm ord}(A_{\alpha})={\rm ord}(\mathcal{Y}_{\alpha})=k_{\alpha}\in \mathbb{N}$ by Lemma \ref{property of order}-(3).

\medskip

\noindent {\bf Claim 1}. ${\rm ord}(\mathcal{Y}_{M, \alpha})\geq {\rm ord}(\mathcal{Y}_{\alpha})$.
\begin{proof}[Proof of Claim 1]
Note that for each $A\in \mathcal{Y}_{\alpha}$, ${\rm ord}(A)\geq {\rm ord}(\mathcal{Y}_{\alpha})$. For each $A\in 2^Y$, if $(x,A)\in W_{\alpha}$, then $A\subset J_{\alpha}[x]$.

Now by Lemma \ref{property of order}~(3), we can choose $B\in \mathcal{Y}_{M, \alpha}$ such that
\[ {\rm ord}(\mathcal{Y}_{M, \alpha})={\rm ord}(B)\ \ \text{and}\ \ B=J_{\alpha}[x_0] \]
for some $x_0\in M$.  Take $(x_0,A)\in W_{\alpha}$. Then $A\in \mathcal{Y}_{\alpha}$, $A\subset B$, and hence ${\rm ord}(A)\leq {\rm ord}(B)$. Thus,
${\rm ord}(\mathcal{Y}_{M, \alpha})={\rm ord}(B)\geq {\rm ord}(A)\geq {\rm ord}(\mathcal{Y}_{\alpha})$.
\end{proof}

\noindent {\bf Claim 2}. For any $\alpha\neq \beta\in \Gamma$, $\Y_\a\perp \Y_{M_i,\b}$ for all $i\in \N$.

\begin{proof}[Proof of Claim 2] Consider $\alpha\neq \beta\in \Gamma$.
Let $M$ be a minimal set  in $(X_{\beta},T_{\beta})$ with $X_{\alpha}\perp M$. By Lemma \ref{disjoint for almost 1-1} $W_{M, \beta}\perp W_{\alpha}$ because both $\pi_1^\beta: (W_{M, \beta},T_\beta\times S)\rightarrow (M,T_\beta)$ and $\pi_1^\alpha: (W_{\alpha},T_{\alpha})\rightarrow (X_{\alpha},T_{\alpha})$ are almost one-to-one extensions.  Consequently, $\mathcal{Y}_{\alpha}\perp\mathcal{Y}_{M,\beta}$ by Lemma \ref{easy}.
\end{proof}

\noindent {\bf Claim 3}. For any $\alpha\in \Gamma$, if there are minimal sets $\{M_i\}_{i\in\mathbb{N}}$ in $X_{\a}$ with $\overline{\bigcup_{i=1}^{\infty}M_i}=X_{\a}$, then $\Y_\a\subseteq \overline{\bigcup_{i=1}^\infty\Y_{M_i,\a}}$.

\begin{proof}[Proof of Claim 3]
Assume that there are minimal sets $\{M_i\}_{i\in\mathbb{N}}$ in $X_{\a}$ with $\overline{\bigcup_{i=1}^{\infty}M_i}=X_{\a}$.
To show $\Y_\a\subseteq \overline{\bigcup_{i=1}^\infty\Y_{M_i,\a}}$, it suffices to show that the transitive point $A_\a=J_\a[x_\a]\in \overline{\bigcup_{i=1}^\infty\Y_{M_i,\a}}$.

Since $\overline{\bigcup_{i=1}^{\infty}M_i}=X_{\a}$ and $\Omega_{M_i,\a}$ is dense in $M_i$ for each $i\in \N$, we have $\overline{\bigcup_{i=1}^{\infty}\Omega_{M_i,\a}}=X_{\a}$. Thus there are $i_m\in \mathbb{N}$ and $x_{i_m}\in \Omega_{M_{i_m},\a}$ for each $m\in \mathbb{N}$  such that $\lim_{m\to \infty} x_{i_m}=x_{\a}$. Since $x_\a\in \Omega_\a$, one has $J_{\a}[x_{i_m}]\rightarrow J_{\a}[x_{\a}]=A_{\a}$ in $2^{Y}$. Note that $J_\a[x_{i_m}]\in \Y_{M_i,\a}$ as $x_{i_m}\in \Omega_{M_{i_m},\a}$. It follows that $A_\a=J_\a[x_\a]\in \overline{\bigcup_{i=1}^\infty\Y_{M_i,\a}}$.
\end{proof}

To sum up, we have a collection of non-trivial transitive subsystems $\{\Y_\a\}_{\a\in \Gamma}$ of $(2^Y, S)$ such that for each $\a\neq \b\in \Gamma$, there exists a family of quasifactors $\{\Y_i\}_{i\in \N}$ such that $\Y_\b\subseteq \overline{\bigcup_{i=1}^\infty \Y_i}$ in $2^Y$, ${\rm ord}(\mathcal{Y}_i)\geq {\rm ord}(\mathcal{Y}_\b)$ and $\Y_\a\perp \Y_i$ for each $i\in\mathbb{N}$.
Then by Theorem \ref{lem-tran-2} $\Gamma$ is at most countable. The whole proof is complete.
\end{proof}



\section{Proof of Theorem \ref{thm-tran}}\label{p-thmA}
In this section, we  give the proof of Theorem \ref{thm-tran} and present a simpler proof when $(X,T)$ is transitive. To do this,  we need the following lemma which is a direct corollary of Theorem \ref{lem-tran-2}.

\begin{lem}\label{yyyy}
Let $(X,T)$ be a non-trivial minimal system and $(\mathcal{W},T)$ be a non-trivial transitive subsystem of $(2^X,T)$.
If $\{(\mathcal{X}_i,T))_{i\in\mathbb{Z}}$ is a sequence of quasifactors such that $\mathcal{W}\subset\overline{\bigcup_{i=1}^\infty \mathcal{X}_i}$ and ${\rm ord}(\mathcal{X}_i)\geq {\rm ord}(\mathcal{W})$ for each $i\in\mathbb{N}$, then
$\mathcal{W}\not\perp \mathcal{X}_{i_0}$ for some $i_0\in\mathbb{N}$.
\end{lem}

\begin{proof}
Otherwise, suppose that $\mathcal{W}\perp\mathcal{X}_i$ for each $i\in \mathbb{N}$.
Let $\Lambda$ be an uncountable index set. For each $\alpha\in \Lambda$, let $\X_{\alpha}={\mathcal W}$. This collection of non-trivial transitive subsystems $\{(\X_\alpha,T)\}_{\alpha\in \Lambda}$ of $(2^X,T)$ satisfies the conditions in Theorem \ref{lem-tran-2}. It follows that 
$\Lambda$ is countable, a contradiction!
\end{proof}


\subsection{Proof of Theorem \ref{thm-tran}}
Now we are going to prove Theorem \ref{thm-tran}.
\begin{proof}[Proof of Theorem \ref{thm-tran}] Let $(X,T)$ be a system. Firstly, assume that $X\perp \mathcal{M}$.  Hence, by Theorem~ \ref{dense-mi}, there are minimal subsystems  $\{M_i\}_{i\in\mathbb{N}}$ of $(X,T)$ such that each $M_i$ is disjoint from $X$ and $\bigcup_{i\in \mathbb{N}}M_i$ is dense in $X$.

\medskip

Conversely, assume that there are minimal subsystems  $\{M_i\}_{i\in\mathbb{N}}$ of $(X,T)$ such that each $M_i$ is disjoint from $X$ and $\bigcup_{i\in \mathbb{N}}M_i$ is dense in $X$. We are going to show that $X\perp \mathcal{M}$.
If this is not true, then there exists a minimal system $(Y,S)$ such that $X\not \perp Y$. Clearly, $Y$ is non-trivial.  Let $J$ be a nontrivial joining of $(X,T)$ and $(Y,S)$, i.e., $J\neq X\times Y$.

Let $X_J$ be the set of all continuous points of the u.s.c map  $\varphi_J: x\mapsto J[x]$.
Next, we consider the map $\Phi: X\rightarrow \mathbb{N}\cup \{\infty\}$ defined  by $\Phi(x)={\rm ord}(J[x])$ for any $x\in X$. Since $\varphi_J$ and ${\rm ord}$ are u.s.c., it follows that $\Phi={\rm ord}\circ \varphi$ is u.s.c. Let $X_\Phi$ be the set of all continuous points the u.s.c map  $\Phi$. Then $X_\Phi$ and $X_J$ are both dense $G_\delta$ sets of $X$ by Lemma \ref{continuous-point-K68}.

For each $i\in \mathbb{N}$, let $E_i:=\left\{x\in X_0: (\overline{\O}(x,T),T)\perp (M_i,T)\right\}$, where $X_0$ is the set of continuous points for $x\mapsto \overline{\O}(x,T)$. Since $X\perp M_i$, by Theorem \ref{dense of disjoint of orbits}-(1), one has $E_i$ is a dense $G_\delta$ subset of $X$.
Now let $E:=\bigcap_{i=1}^\infty E_i$, then $E$ is also a dense $G_\delta$ subset of $X$. Thus $ X_0\cap E\cap X_\Phi\cap X_J$ is a dense $G_\delta$ set of $X$.


Since $J\neq X\times Y$, there is a nonempty open set $U\subset X$ such that
\[ J[x]\neq Y, \ \ \text{ for any } x\in U.\]
Hence, there is some $x^{*}\in U\cap X_0\cap E\cap  X_\Phi\cap X_J$ with $ J[x^{*}]\neq Y$.
Thus, by Lemma \ref{property of order}-(2), $\Phi(x^*)={\rm ord}(J[x^*])<+\infty$. Since $\Phi:X\rightarrow \mathbb{N}\bigcup \{\infty\}$ is continuous at $x^*$, there exists an open neighborhood $V$ of $x^*$ such that $V\subseteq U$ and $\Phi(x)=\Phi(x^*)<+\infty$ for any $x\in V$.

Since $\overline{\bigcup_{i=1}^{\infty}M_{i}}=X$, one has $V\cap (\bigcup_{i=1}^{\infty}M_{i})$ is dense in $V$.  Thus, there exist a sequence $\{n_i\}_{i\in \mathbb{N}}$ of $\mathbb{N}$ and $\{x_i\}_{i\in \mathbb{N}}$ of $X$ such that
\begin{itemize}
\item  $x_i\in V\cap M_{n_i}$ for each $i\in \mathbb{N}$.
\item $\lim_{i\to \infty}x_i=x^*$.
\end{itemize}
Since $x^*\in X_0$ and $\lim_{i\to \infty}x_i=x^*$, one has
$$\lim_{i\to \infty} M_{n_i}=\lim_{i\to \infty}\overline{\O}(x_i,T)=\overline{\O}(x^*,T) \text{ in }2^X.$$

Now we use Corollary \ref{almos 1-1 induced by joining} to induce systems in $(2^Y,S)$.
Let $\Omega_{*}$ be the set of continuous points of the map $\overline{\O}(x^{*}, T)\rightarrow 2^{Y}: x\mapsto J[x]$ and let $\Omega_i$ be the set of continuous points of the map $M_{n_i}\rightarrow 2^{Y}: x\mapsto J[x]$ for each $i\in \mathbb{N}$. Let
\[ W_{*}=\overline{\{ (x, J[x])\in \overline{\O}(x^{*}, T)\times 2^{Y}: x\in\Omega_{*}\}},\]
\[ W_i=\overline{\{ (x, J[x])\in M_{n_i}\times 2^{Y}: x\in\Omega_i\}}, \forall i\in \N,\]
and \[ \mathcal{W}=\pi_2(W_{*})\ \ \text{and} \ \ \mathcal{X}_{i}=\pi_2(W_i).\]
Then the first coordinate projection $\pi_1:(W_*,T\times S)\rightarrow (\overline{\O}(x^{*}, T),T)$ and $\pi_1: (W_i,T\times S)\rightarrow (M_{n_i},T)$ are both  almost one-to-one extensions. By Lemma \ref{almost1-1}, $(W_*, T\times S)$ and $({\mathcal W}, S)$ are transitive as $(\overline{\O}(x^{*}, T),T)$ is transitive; for each $i\in \mathbb{N}$, $(W_i, T\times S)$ and $\X_i$ are minimal as $(M_{n_i},T)$ is minimal. Note that $Y\neq J[x^*]\in \mathcal{W}$ as $x^*\in U$. Hence, $(\mathcal{W},S)$ is a non-trivial transitive system
because the unique fixed point in $(2^Y,S)$ is $Y$.

Since $x^*\in E$, we have $\overline{\O}(x^{*}, T)\perp M_{n_i}$, and hence $W_*\perp W_i$  by Lemma \ref{disjoint for almost 1-1} for each $i\in \mathbb{N}$, and $\mathcal{W}\perp \mathcal{X}_i$ for each $i\in \mathbb{N}$.

\medskip

We claim that

\medskip

\noindent{\bf Claim.} For each $i\in \mathbb{N}$, ${\rm ord}(\mathcal{X}_{i})= {\rm ord}(\mathcal{W})$, and
$\mathcal{W}\subset \overline{\bigcup_{i=1}^{\infty}  \mathcal{X}_{i}}$ in $2^Y$.

\begin{proof}[Proof of Claim]
Since $x^*\in X_J\cap \overline{\O}(x^*,T)\subset \Omega_*$, we have $\overline{\O}((x^*,J[x^*]),T\times S)=W_*$, and thus
$\overline{\O}(J[x^*],S)=\mathcal{W}$. Hence, by Lemma \ref{property of order}-(3), ${\rm ord}(\mathcal{W})={\rm ord}(J[x^*])=\Phi(x^*)<\infty$.
For each $i\in \mathbb{N}$, since $x_{n_i}\in V\cap M_{n_i}$  and $\Omega_i$ is a dense $G_\delta$ set of $M_{n_i}$, we may take $z_i\in V\cap M_{n_i}\cap \Omega_i$ such that
$\rho(z_i,x_{n_i})<\frac{1}{2^i}$, where $\rho$ is the metric on $X$. Then $J[z_i]\in  \mathcal{X}_{i}$, and thus $${\rm ord}(\mathcal{X}_{i})={\rm ord}(J[z_i])=\Phi(z_i)=\Phi(x^*)$$
as $z_i\in V$. This implies that ${\rm ord}(\mathcal{X}_{i})={\rm ord}(\mathcal{W})$.

Next, since $\overline{\O}(x^*,T)\subset X_J\cap \overline{\O}(x^*,T)\subset \Omega_*$  and  $\O(x^*,T)$ is dense in $\Omega_*$, we have
$$\mathcal{W}=\overline{\{J[x]:x\in \Omega_*\}}= \overline{\{J[T^nx^*]:n\in \mathbb{Z}\}}=\overline{\{S^n(J[x^*]):n\in \mathbb{Z}\}}.$$
Note that $S(\overline{\bigcup_{i=1}^{\infty}  \mathcal{X}_{i}})=\overline{\bigcup_{i=1}^{\infty}  \mathcal{X}_{i}}$. Hence, to show $\mathcal{W}\subset \overline{\bigcup_{i=1}^{\infty}  \mathcal{X}_{i}}$, it suffices to show that $J[x^*]\in \overline{\bigcup_{i=1}^{\infty}  \mathcal{X}_{i}}$.
Since $z_i\in \Omega_i$, we have $J[z_i]\in \mathcal{X}_i$ for each $i\in \mathbb{N}$. Note that
$$\lim_{i\to\infty}z_i=\lim_{i\to\infty}x_{n_i}=x^*\in X_J.$$
Thus, $J[x^*]=\lim_{i\to\infty}J[z_i]\in \overline{\bigcup_{i=1}^{\infty}  \mathcal{X}_{i}}$.
The proof of Claim is complete.
\end{proof}

Summing up, we have shown that $\mathcal{W}$  and $\mathcal{X}_i$'s satisfy the conditions in Lemma \ref{yyyy} by the above Claim. But then we get a contradiction with the above fact that $ \mathcal{W}\perp  \mathcal{X}_{i}$ for each $i\in \mathbb{N}$, because Lemma \ref{yyyy} tells us that $ \mathcal{W}\not\perp  \mathcal{X}_{i_0}$ for some $i_0\in \mathbb{N}$. Therefore, we have showed $X\perp Y$ and hence $X\perp\mathcal{M}$. The whole proof is complete.
\end{proof}

\subsection{A simpler proof of Theorem \ref{thm-tran} for a transitive system}\label{appen-B}
When $(X,T)$ is transitive, we have a much simpler proof which does not involve quasifactors. We start from the following lemma.
\begin{lem}\label{lift of min pt}
Let $(X,T)$ and $(Y,S)$ be two minimal systems.
\begin{enumerate}
\item[(1)] For each $y\in Y$, there is some point $x\in X$ such that $(x,y)$ is a minimal point in $(X\times Y, T\times S)$.
\item[(2)] If $(x,y)$ is a minimal point of $(X\times Y,T\times S)$, then for any $m,n\in\mathbb{Z}$, $(T^{n}x, S^{m}y)$ is also a minimal point of $(X\times Y,T\times S)$.
\end{enumerate}
\end{lem}
 \begin{proof}
 (1) It is east to see that $p_2: X\times Y\rightarrow Y, (x,y)\mapsto y$  is a factor map between $(X\times Y, T\times S)$ and $(Y,S)$.
 Hence (1) follows from this fact.

 (2) For each $m,n\in\mathbb{Z}$, $T^{n}\times S^{m}$ commutes with $T\times S$. Thus if $(x,y)$ is a minimal point of $(X\times Y,T\times S)$, then $(T^{n}x, S^{m}y)$ is also a minimal point of $(X\times Y,T\times S)$.
 \end{proof}

 Now we are ready to show the following theorem. 
 \begin{thm} Let $(X,T)$ be a transitive system. Then the following statements are equivalent.
 \begin{enumerate}
 \item $X\perp \mathcal{M}.$
 \item There exist minimal subsets $\{M_i\}_{i \in \mathbb{N}}$ in $X$  with $\overline{\cup_{i \in \mathbb{N}} M_i} = X$  such that each $M_i$ is disjoint from $X$.
 \item There exist minimal subsets $\{M_i\}_{i \in \mathbb{N}}$ in $X$  with $\overline{\cup_{i\in \mathbb{N}} M_i} = X$  such that $(x,y)$
 is proximal for any transitive point $x$ of $X$ and any $y\in \cup_{i \in \mathbb{N}} M_i.$
 \item $(X,T)$ has dense minimal points and $(x,y)$ is proximal for any transitive point $x$ of $X$ and any minimal point $y$ of $X$.
 \end{enumerate}
 \end{thm}
\begin{proof}
(1) $\Rightarrow$ (2). It follows from Theorem \ref{dense-mi}.

(2) $\Rightarrow$ (3). Let $M_i$ be the minimal set in (2) for any $i\in\N$, $x$ be a transitive point and $y\in M_i$ for some $i\in \N$.
It is clear that $\overline{\O}((x,y),T\times T)$ is a joining of $X$ and $M_i$ and hence $\overline{\O}((x,y),T\times T)=X\times M_i$, as $X\perp M_i$. Since $\overline{\O}((x,y),T\times T)=X\times M_i\supset M_i\times M_i\supset \{(t,t):t\in M_i\}$, we conclude that $(x,y)$ is proximal.

(3) $\Rightarrow$ (1) Now suppose that there exist minimal subsets $\{M_i\}_{i \in \mathbb{N}}$ in $X$  with $\overline{\cup_{i\in \mathbb{N}} M_i} = X$  such that $(x,y')$  is proximal for any transitive point $x$ of $X$ and any $y'\in \cup_{i \in \mathbb{N}} M_i.$
To show $X\perp\mathcal{M}$, we fix a minimal system $(Y,S)$ and aim to show  $X\perp Y$.

Let $J$ be a joining of $(X,T)$ and $(Y,S)$. Fix a transitive point $x\in X$. Then we can find $y\in Y$ with $(x,y)\in J$.
By Lemma \ref{lift of min pt}~(1), there is  $x_i\in M_i$ such that $(x_i,y)$ is a minimal point of $(M_i\times Y,T\times S)$ and hence of $(X\times Y,T\times S)$, for each $i\in \mathbb{N}$. Further, by Lemma \ref{lift of min pt}~(2), $(T^{n}x_i,y)$ is also a minimal point in $(M_i\times Y,T\times S)$ for each $n\in\mathbb{Z}$. Set
\[ D=\{ T^{n}x_i: i\in\mathbb{N}, n\in\mathbb{Z}\},\]
which is clearly dense in $X$.

Now for each $z\in D$, it is clear that $(x,z)$ is proximal. We fix $z\in D$ and $\epsilon>0$.  Then the set
$C_1=\{n\in \mathbb{Z}: \rho(T^{n}x, T^{n}z)<\epsilon\}$
contains arbitrarily long intervals of $\Z$.
On the other hand, since $(z,y)$ is a minimal point of $(X\times Y,T\times S)$, the set
$C_2=\{n\in \mathbb{Z}: \rho(T^{n}z,z)<\epsilon, \rho(S^{n}y,y)<\epsilon\}$
is syndetic, i.e. of bounded gaps. Thus $C_1\cap C_2\neq\emptyset$. In other words, we can find $n\in \mathbb{Z}$ such that
\[ \rho(T^{n}x, T^{n}z)<\epsilon,\ \ \rho(T^{n}z,z)<\epsilon,\ \  \rho(S^{n}y,y)<\epsilon.\]
Hence $\rho(T^{n}x, z)<2\epsilon$ and $\rho(S^{n}y,y)<\epsilon.$
Since $\epsilon>0$ is arbitrary, we get $ (z,y)\in \overline{\O}((x,y),T\times S)$. As $z$ is also chosen arbitrarily in $D$ and $D$ is dense in $X$, we get that $X\times \{y\}\subset \overline{\O}((x,y),T\times S)\subset J.$
Since $(Y,S)$ is minimal, we have $J=X\times Y$.  This shows $X\perp Y$.

(1) $\Rightarrow$ (4). By Theorem \ref{dense-mi} $(X,T)$ has dense minimal points. Now let $y\in X$ be a minimal point. Then $X\perp \overline{\O}(y,T)$. So the  proof of (2) $\Rightarrow$ (3) implies that $(x,y)$ is proximal.

(4) $\Rightarrow$ (1). It is clear that (4) $\Rightarrow$ (3) and hence (4) $\Rightarrow$ (1).
\end{proof}
By the above theorem, if a transitive system is disjoint from any minimal system, then any transitive point is proximal to any minimal point in the system. We remark that generally a transitive point of a transitive system $(X,T)$ is only proximal to some minimal points in $X$, see for example \cite{Au88}.


\section{The proofs of Theorem \ref{thm general} and Theorem \ref{thm-ss}}\label{sec-gen}

We prove Theorem \ref{thm general}  and Theorem \ref{thm-ss} in this section.

\subsection{Some preparations}

First, we need the following topological ergodic decomposition given by Akin and Glasner.

\begin{thm}[Topological ergodic decomposition, {\cite[Theorem 1]{AG98}}]\label{TED}
Let $(X,T)$ be a topological dynamical system. Then there is an almost one-to-one extension $\theta: (X^{*},T^*)\rightarrow (X,T)$ which has a factor $\pi_*: (X^{*}, T^*) \rightarrow (Z, \id)$ satisfying
\begin{enumerate}
\item[(1)] $\pi_*$ is open and
\item[(2)] there is a dense $G_{\delta}$ subset $\Omega$ of $X^*$ such that for each $x\in \Omega$, $\pi_*^{-1}(\pi_*(x))=\overline{\O}(x,T^*)$.
\end{enumerate}
\end{thm}

It is clear that the factor $(Z, \id)$ is an equicontinuous factor of $(X^{*}, T^{*})$ in Theorem \ref{TED}. But in general it is not the maximal equicontinuous factor. If $(X,T)\perp \mathcal{M}$ then we will show it is the case. Now we show that if $(X,T)$ has DDMS-property, then the maximal equicontinuous factor $X_{eq}$ of $(X,T)$ consists solely of fixed points.

\begin{thm}\label{ge-sufficient-condition} Let $(X,T)$ be a system. If $(X,T)$ has DDMS-property, then $$\overline{\O}(x,T)\times \overline{\O}(x,T)\subset Q(X)$$ for each $x\in X$.
Consequently, the maximal equicontinuous factor $X_{eq}$ of $(X,T)$ consists solely of fixed points.
\end{thm}

\begin{proof}
Assume that $(X,T)$ has DDMS-property, i.e. there are minimal subsystems  $\{X_i\}_{i\in\mathbb{N}}$ of $(X,T)$ such that each $X_i$ is disjoint from $X$ and $\bigcup_{i\in \mathbb{N}}X_i$ is dense in $X$.

First, we show that $X_i\times X_i\subset Q(X)$ for each $i\in \mathbb{N}$. When $X_i$ is singleton, it is clear that $X_i\times X_i\subset Q(X)$.
Now, let us consider a non-trivial minimal subsystem  $(X_i,T)$ of $(X,T)$ with $X_i\perp X$. Since $X_i\perp X$, by Theorem \ref{hy-general} there are maximal transitive subsystems $(Y_j,T)$ in $(X,T)$ such that $(Y_j,T)\bot (X_i,T), j\in\mathbb{N}$ and $\bigcup_{j\in\mathbb{N}} Y_j$ is dense in $X$.

Let $x,y\in X_i$. We need to show that $(x,y)\in Q(X)$. Since $\overline{\bigcup_{j=1}^\infty Y_j}=X$ and the set of transitive points in a transitive system is a residual set, there are sequence $\{j_m\}_{m=1}^\infty\subseteq \N$ and transitive points $y_m\in Y_{j_m}, m\in \mathbb{N}$ such that $\lim_{m\to\infty}y_{j_m}=x$. Without loss of generality assume that $\lim_{m\to \infty}Y_{j_m}=M$ in $2^X$, then $M\in 2^X$ is $T$-invariant and $X_i\subset M$.

Let $U\times V$ be an open neighborhood of $(x,y)$. Since $X_i\subset M=\lim_{m\to \infty}Y_{j_m}$, there is some $m\in \mathbb{N}$ such that $Y_{j_m}\cap V\neq\emptyset$. Take $y'\in Y_{j_m}\cap V$. Then  $(x,y')\in (U\times V)\cap X_i \times Y_{j_m}, (y, y')\in (V\times V)\cap X_i\times Y_{j_m}$. This implies that $(U\times V)\cap X_i \times Y_{j_m}, (V\times V)\cap X_i\times Y_{j_m}$ are two non-empty open subset of $X_i\times Y_{j_m}$. Since $X_i \bot Y_{j_m}$, $(X_i\times Y_{j_m}, T\times T)$ is transitive.
Hence there is a transitive point $(z_1,z_2)$ of the system $(X_i\times Y_{j_m},T\times T)$ and $n \in \mathbb{Z}$ such that  $(z_1,z_2)\in (U\times V)\cap X_i \times Y_{j_m}$ and  $(T^nz_1,T^nz_2)\in (V\times V)\cap X_i \times Y_{j_m}$. This implies that $(x,y)\in Q(X)$. By the arbitrariness of $x,y$, we conclude $X_i\times X_i\subset Q(X)$.

\medskip

Now we show that $\overline{\O}(x,T)\times \overline{\O}(x,T)\subset Q(X)$ for all $x\in X$. Since $\bigcup_{i\in \mathbb{N}}X_i$ is dense in $X$, there are sequence $\{i_\ell\}_{\ell=1}^\infty\subseteq\N$ and $\{x_\ell\}_{\ell=1}^\infty \subset X$ such that $x_\ell\in X_{i_\ell}$ for each $\ell\in \mathbb{N}$, and  $\lim_{\ell\to \infty}x_\ell=x$. Without loss of generality suppose that $\lim_{\ell\to\infty}X_{i_\ell}\to Y'$ in $2^X$, then $Y'\in 2^X$ is $T$-invariant and $x\in Y'$. Thus $\overline{\O}(x,T)\subset Y'$, and hence $\overline{\O}(x,T)\times \overline{\O}(x,T)\subset Y'\times Y'$. Note that for each minimal subsystem $(X_{i_\ell},T)$, $X_{i_\ell}\times X_{i_\ell}\subset Q(X)$ by the previous discussion. By letting $\ell \to \infty$, we have $Y'\times Y'\subset Q(X)$ as $Q(X)$ is closed. In particular, we have that $\overline{\O}(x,T)\times \overline{\O}(x,T)\subset Q(X)$.

\medskip
Assume that $X_{eq}$ is the maximal equicontinuous factor of $X$. For each minimal subsystem $Y$ of $X_{eq}$ there is a minimal subsystem $Y'$ of $X$
with $\pi(Y')=Y$, where $\pi:X\ra X_{eq}$ is the factor map. Since $Y'\times Y'\subset Q(X)$ we conclude that $Y$ is a fixed point. The proof is complete.
\end{proof}

By Theorem \ref{dense-mi}, $(X,T)\bot\mathcal{M}$ implies that $(X,T)$ has DDMS-property. Thus we have the following corollary.

\begin{cor}\label{dis-sufficient-condition}
Let $(X,T)$ be a system. If $(X,T)\bot\mathcal{M}$, then
$\overline{\O}(x,T)\times \overline{\O}(x,T)\subset Q(X)$
for each $x\in X$.
Consequently, the maximal equicontinuous factor $X_{eq}$ of $X$ is consisting of fixed points.
\end{cor}


\subsection{Proof of Theorem \ref{thm general}}
Now we are ready to give the proof of Theorem \ref{thm general}.

\begin{proof}[Proof of Theorem \ref{thm general}]
$(\Longrightarrow)$  Assume $X\perp\mathcal{M}$. We consider an almost one-to-one extension $(X^{*},T^*)$ of $(X,T)$ as in Theorem \ref{TED}. Consequently, $X^*\perp\mathcal{M}$ by Lemma \ref{disjoint for almost 1-1}.

Let $\pi_*:(X^{*}, T^*)\rightarrow (Z, \id)$ be the factor map given in Theorem \ref{TED} and $\Omega$ be a dense $G_{\delta}$ subset of $X^*$ such that for each $x\in \Omega$, $\pi_*^{-1}(\pi_*(x))=\overline{\O}(x,T^*)$. We first show that $Z=X^{*}_{eq}$ is the maximal equicontinuous factor of $X^{*}$. Clearly, $Z$ is a factor of $X_{eq}^{*}$ since $Z$ is equicontinuous.
Let $\phi: (X^*_{eq}, T^*_{eq})\rightarrow (Z,\id)$ be the factor map. By Corollary \ref{dis-sufficient-condition}, $$\overline{\O}(x,T^*)\times \overline{\O}(x,T^*)\subset Q(X^*)$$ for each $x\in X^*$, and the maximal equicontinuous factor $X_{eq}^*$ is consisting of fixed points. Since for each $x\in \Omega$, $\pi_*^{-1}(\pi_*(x))=\overline{\O}(x,T^*)$, we have $\phi$
is almost one-to-one. But as a factor of open factor map $\pi_*$, $\phi$ is also open. Thus $\phi$ is an isomorphism and $X^*_{eq}=Z$. In particular, $\pi: X^{*}\rightarrow X^{*}_{eq}$ is a topological decomposition.

Since $X^*\perp \mathcal{M}$, $(X^*,T^*)$ has DDMS-property by Theorem \ref{thm-tran}. Hence there are minimal sets $\{M_i\}_{i\in\mathbb{N}}$ in $X^{*}$ such that $M_i\perp X^*$ and $\overline{\bigcup_{i=1}^{\infty}M_i}=X^{*}$. By Corollary \ref{many pt for disjoint} , for each $i\in\mathbb{N}$, the set $\{x\in X^{*}: \overline{\O}(x, T^*)\perp M_i \}$ is a residual set of $X^*$ because $M_i\perp (X^*,T^*)$. Let $Rec(X^{*})$ be the set of all recurrent points of $(X^{*},T^*)$, which is a dense $G_{\delta}$ subset of $X^*$. Meantime, the map $\overline{\O}: X^*\ra 2^{X^*}$, $x\mapsto \overline{\O}(x,T^*)$ is l.s.c, and the set
$\Omega_0$ of the set of continuous point of $\overline{\O}$ is a dense $G_\delta$ subset of $X^*$.
Thus, the set
\[\Omega_1=\Omega_0\cap Rec(X^{*})\cap \bigcap_{i\in \mathbb{N}}\{ x\in X^*: \overline{\O}(x, T^*)\perp M_i\}\]
is a residual set of $X^{*}$.
Let $\Omega_2$ be the set of continuous points of the map
\[ \Phi: X^{*}\times X^{*}\rightarrow 2^{X^{*}\times X^{*}},  (x,y)\mapsto \overline{\O}((x,y),T^*\times T^*).\]
Then $\Omega_2$ is a dense $G_{\delta}$ set of $X^{*}\times X^{*}$. Further, set
\[ \Omega_3=\Omega_2\cap(\Omega_1\times \Omega_1),\]
which is a residual set of $X^{*}\times X^{*}$.  Finally, we set
\[Z^{*}=\pi\times \pi(\Omega_3\cap (\Omega\times \Omega)),\]
which is a residual set of $X^{*}_{eq}\times X^{*}_{eq}$ by Theorem \ref{veech} since $\pi\times \pi$ is an open map.


Let $(z_1,z_2)\in Z^{*}$. As $z_1,z_2\in \pi(\Omega_1)$, both $\pi^{-1}(z_1)$ and $\pi^{-1}(z_2)$ are transitive and are disjoint from each $M_i$, i.e. $\pi^{-1}(z_i)\perp M_j$ for any $i=1,2$ and $j\in\mathbb{N}$.

It remains to show that $(\pi^{-1}(z_1)\times \pi^{-1}(z_2),T^*\times T^*)$ is a transitive subsystem for each $(z_1,z_2)\in Z^{*}$.
This is followed by fact $\pi^{-1}(z_i)\perp M_j$ for any $i=1,2$ and $j\in\mathbb{N}$ and the continuity of $\Phi$ on $\Omega_2$.

To be precise, fix $(z_1,z_2)\in Z^{*}$ and take $(y_1,y_2)\in \Omega_3\cap (\Omega\times \Omega)\cap(\pi^{-1}(z_1)\times\pi^{-1}(z_2))$. Since $y_1,y_2\in \Omega$, one has $\pi^{-1}(z_1)=\overline{\O}(y_1,T^*)$ and  $\pi^{-1}(z_2)=\overline{\O}(y_2,T^*)$. Since $\overline{\bigcup_{i=1}^{\infty}M_i}=X^{*}$, we may take sequence $\{n_i\}_{i\in \mathbb{N}}$ of $\mathbb{N}$ and $x_{n_i}\in M_{n_i}, i\in \mathbb{N}$ such that $\lim_{i\to \infty} x_{n_i}=y_2$. In particular, $(y_1,x_{n_i})\to (y_1,y_2)$ as $i\to\infty$. Since $\Phi$ is continuous at $(y_1,y_2)$, one has
\[ \overline{\O}((y_1,y_2),T^*\times T^*)=\lim_{n_i\to \infty} \overline{\O}((y_1,x_{n_i}),T^*\times T^*). \]
But $\overline{\O}((y_1,x_{n_i}),T^*\times T^*)=\overline{\O}(y_1,T^*)\times M_{n_i}$ as $\overline{\O}(y_1,T^*)\perp M_{n_i}$. Clearly, $M_{n_i}\to \overline{\O}(y_2, T^*)$ in $2^{X^*}$ as $\lim_{i\to \infty}x_{n_i}=y_2$ and $y_2\in \Omega_0$. Thus
\[ \overline{\O}((y_1,y_2),T^*\times T^*) =\overline{\O}(y_1, T^*)\times\overline{\O}(y_2, T^*)=\pi^{-1}(z_1)\times \pi^{-1}(z_2). \]
In particular, $(\pi^{-1}(z_1)\times \pi^{-1}(z_2),T^*\times T^*)$ is a transitive system.

\medskip
$(\Longleftarrow)$ Assume that $(X,T)$ has an almost one-to-one extension $(X^{*},T^*)$ such that the maximal equicontinous factor map $\pi: X^{*}\rightarrow X^{*}_{eq}$ is a topological decomposition, and there are countably many minimal sets $\{M_i\}_{i\in\mathbb{N}}$ in $(X^{*},T^{*})$ with $\overline{\bigcup_{i=1}^{\infty}M_i}=X^{*}$ and a dense $G_{\delta}$ set $Z^{*}$ of $X_{eq}^{*}\times X_{eq}^{*}$ such that for each $(z_1,z_2)\in X^{*}$,
\begin{itemize}
\item[(a)]  $(\pi^{-1}(z_1)\times \pi^{-1}(z_2), T^*\times T^*)$ is a transitive system;
\item[(b)] $\pi^{-1}(z_i)\perp M_j$ for any $i=1,2$ and $j\in\mathbb{N}$.
\end{itemize}

Let $\pi_1: X_{eq}^*\times X_{eq}^*\rightarrow X_{eq}^*$ denote the projection onto the first coordinate. Since $\pi_1$ is open, $\pi_1(Z^*)$ is a residual set of $X_{eq}^*$ by Theorem \ref{veech}. Since $\pi$ is open, $\pi^{-1}(\pi_1(Z^{*}))$ is a residual set of $X^*$ by Lemma \ref{residual under semiopen}.
Next, for each $z\in \pi_1(Z^*)$, by the above (a) and (b), $(\pi^{-1}(z),T^*)$ is transitive and $\pi^{-1}(z)\perp M_j$ for each $j\in \mathbb{N}$.
Now for each $j\in \mathbb{N}$, since $\pi^{-1}(z)\perp M_j$ for each $z\in \pi_1(Z^*)$ and $\pi^{-1}(\pi_1(Z^*))=\bigcup_{z\in \pi_1(Z^*)}\pi^{-1}(z)$ is a residual set of $X^*$, we have $\{x\in X^*: \overline{\O}(x,T^*)\perp M_j\}$ is a residual subset of $X^*$. It follows that $M_j\perp X^*$ by Theorem \ref{dense of disjoint of orbits}-(2). Therefore, by Theorem \ref{thm-tran}, $X^*\perp \mathcal{M}$, and thus $X\perp \mathcal{M}$ by Lemma \ref{easy}-(1). The proof is complete.
\end{proof}

\subsection{Proof of Theorem \ref{thm-ss}}
Now we are going to prove Theorem \ref{thm-ss}. Let  $(X,T)$ be a semi-simple system. We are going to show that $X\bot\mathcal{M}$ if and only if  the set
\[\Delta_{\perp}(X):= \left\{(x_1,x_2)\in X\times X: (\overline{\O}(x_1,T),T)\perp (\overline{\O}(x_2,T),T)\right\}\]
is residual in $X\times X$.

In Corollary \ref{dense G of disjoint of orbits}, we have show one direction: if  $X\bot\mathcal{M}$, then $\Delta_{\perp}(X)$ is residual in $X\times X$. Now we show another direction.

Assume that $\Delta_{\perp}(X)$ is residual in $X\times X$. Then there exist dense open subsets $\{V_n\}_{n=1}^\infty$ of $X\times X$ such that
$\Delta_{\perp}(X)\supset \bigcap_{n=1}^\infty V_n$. Firstly, we have

\medskip
\noindent{\bf Claim}. If $z\in X$ is an isolated point in $X$, then $z$ is a fixed point of $(X,T)$.
\begin{proof}[Proof of Claim]
Assume $z$ is an isolated point in $X$. Then $(z,z)$ lies in every dense $G_{\delta}$ set of $X\times X$. In particular, $(z,z)\in\Delta_{\perp}(X)$. But then $(\overline{\O}(z,T),T)\perp (\overline{\O}(z,T),T)$, which implies that $z$ is a fixed point. 
\end{proof}

Now let $(Y,T)$ be a minimal system and we need to show $X\perp Y$.
If this is not true (that is, $X\not\perp Y$), then by Theorem \ref{dense of disjoint of orbits}-(2),  there is a nonempty open set $U\subset X$ such that
\[\overline{\O}(x,T)\not\perp Y,\ \ \forall x\in \overline{U}.\]
By the above Claim, there is no isolated points of $X$ in $U$ because any fixed point is disjoint from $Y$. Hence $\overline{U}$ is a perfect set  of $X$.
Since each $(U\times U)\cap V_n$ is dense open subset of $\overline{U}\times \overline{U}$ for $n\in \mathbb{N}$,
one has that $G:=(U\times U)\cap \bigcap_{n=1}^\infty V_n$ is dense $G_\delta$-set  of $\overline{U}\times \overline{U}$ and $G \subseteq \Delta_{\perp}(X)$.

Applying  Mycielski's Theorem to the perfect set $\overline{U}$ (see for example \cite[Thorem 19.1]{Kechris}), there is a Cantor set $C\subset \overline{U}$ such that
\[ \{ (x_1,x_2)\in C\times C: x_1\neq x_2\}\subset G.\]
Thus for any $x_1\neq x_2\in C$,
\[ \overline{\O}(x_1,T) \perp \overline{\O}(x_2,T).\]
Combining this with the fact that $(\overline{\O}(x,T),T)$ is minimal and $\overline{\O}(x,T)\not\perp Y$ for any $x\in C$,  it follows from Lemma \ref{joining gives QF} that we have a collection $\{\mathcal{X}_{\overline{\O}(x,T)}\}_{x\in C}$ of pairwise disjoint non-trivial quasifactors of $(Y,T)$. But this contradicts with Theorem \ref{thm QF} because $C$ is an uncountable set.
Thus we have that $X\perp Y$. The proof is complete.
\hfill $\square$


\section{Distal systems disjoint from all minimal systems}\label{sec-distal}

In this section, our goal is to characterize the distal systems that are disjoint from all minimal systems. For example, we will show that a distal system is disjoint from any minimal system if and only if the maximal equicontinuous factor of its any almost one-to-one extension consists of fixed points. Also we will prove Theorem \ref{Main-thm-B} in this section.

\subsection{Distal systems}

First we have the following easy lemma about almost one-to-one distal extension.

\begin{lem}\label{easylemmas}
Assume that $(X,T)$ is a distal system and $\pi:(Z,R)\ra (X,T)$ is an almost one-to-one distal extension. Let $Z_0=\left \{z\in Z: \pi^{-1}\pi(z)=\{z\}\right \}$. Then
 for each $z\in Z_0$, we have
\begin{enumerate}
\item $\overline{\O(z,R)}\subset Z_0$ and $\pi^{-1}\left(\overline{\O}(\pi(z),T)\right)=\overline{\O}(z,R)$.
\item $\pi: (\overline{\O}(z,R),R)\rightarrow (\overline{\O}(\pi(z),T),T)$ is an isomorphism.
\end{enumerate}
\end{lem}
\begin{proof} Let $z\in Z_0$. First we show that $\pi^{-1}\left(\overline{\O}(\pi(z),T)\right)=\overline{\O}(z,R)$. Assume the contrary that there is a minimal subsystem $M\subset Z$ such that $M\cap \overline{\O}(z,R)=\emptyset$ and there is $y\in M$ with $\pi(y)\in \overline{\O}(\pi(z),T).$ This implies that
$\pi(M)=\overline{\O}(\pi(z),T)$ and hence $|\pi^{-1} \pi(z)|\ge 2$, a contradiction.

It is clear that $\pi: (\overline{\O}(z,R),R)\rightarrow (\overline{\O}(\pi(z),T),T)$ is almost one-to-one between minimal systems. As $(\overline{\O}(z,R),R)$ is distal, we conclude that $\pi: (\overline{\O}(z,R),R)\rightarrow (\overline{\O}(\pi(z),T),T)$ is a homeomorphism, which implies $\overline{\O(z,R)}\subset Z_0$, since $\pi^{-1}\left(\overline{\O}(\pi(z),T)\right)=\overline{\O}(z,R)$.
\end{proof}


\begin{rem} There are examples where the factor map $\pi:(Z,R)\ra (X,T)$ is almost one-to-one and $(Z,R)$ is distal. For example,  let $Z=[0,1]\times \mathbb{T}$ and define $R(x,y)=(x,x+y)$ for $x\in [0,1]$ and $y\in \mathbb{T}$. Then $(Z,R)$ is distal. Collapsing $\{0\}\times \mathbb{T}$ to a point, we obtain a factor system $(X,T)$ and the corresponding factor map $\pi:(Z,R)\ra (X,T)$. Clearly, $\pi$ is an almost one-to-one distal extension.
\end{rem}

Now we show that $(X,T)\bot \mathcal{M}$ if and only if the maximal equicontinuous factor of any almost one-to-one extension of $X$  consists of fixed points.

\begin{thm}\label{equi-ll}
Let $(X,T)$ be a distal system. The following statements are equivalent:

\begin{enumerate}
\item $(X,T)\bot \mathcal{M}$.

\item The maximal equicontinuous factor of any almost one-to-one extension of $X$  consists of fixed points.

\item The maximal equicontinuous factor of any almost one-to-one distal extension of $X$ consists of fixed points.

\item Every almost one-to-one extension of $(X,T)$ has DDMS-property.

\item Every almost one-to-one distal extension of $(X,T)$ has DDMS-property.
\end{enumerate}
\end{thm}

\begin{proof}
(1) $\Longrightarrow$ (2). Assume that $(X,T)\bot \mathcal{M}$. Let $(X',T')$ be an almost one-to-one extension of $(X,T)$. By Lemma \ref{disjoint for almost 1-1}  $(X',T')\bot \mathcal{M}$. Furthermore, by Corollary \ref{dis-sufficient-condition}, the maximal equicontinuous factor of $(X',T')$ consists of fixed points.

\medskip

(2) $\Longrightarrow$ (3). It is obvious.

\medskip
(3) $\Rightarrow$ (1). Assume that the maximal equicontinuous factor of any almost one-to-one distal extension of $(X,T)$ consists of fixed points.
By Theorem \ref{mini-equi}, to show that $(X,T)\bot \mathcal{M}$ it remains to show that $(X,T)$ is disjoint from any minimal equicontinuous system.

Otherwise assume that there is non-trivial equicontinuous minimal system $(Y,S)$ such that $(X,T)\not\perp (Y,S)$. Thus there is some
joining $J\neq X\times Y$. Since $J\not =X\times Y$, there is $x_0\in X$ such that $J[x_0]\neq Y$, where $J[x_0]=\{y\in Y: (x_0,y)\in J\}$. By Lemma \ref{QF of equi}, $(2^Y,S)$ is equicontinuous with a unique fixed point $\{Y\}$. In particular, $(X\times 2^Y,T\times S)$ is distal.

Let $\Omega\subset X$ be the continuous points of the map $X\rightarrow 2^Y, x\mapsto J[x]$. Let
$$W=\overline{\{(x,J[x]):x\in \Omega\}}\subset X\times 2^Y.$$
$$
\xymatrix{
                & (W, T\times S) \ar[dl]_{\pi_1} \ar[dr]^{\pi_2} \\
 (X,T) & &     (\Y_J,S)        }
$$
By Corollary \ref{almos 1-1 induced by joining}, $\pi_1:(W,T\times S)\ra (X,T), \ (x,A)\mapsto x$ is an almost one-to-one extension of $(X,T)$. Since $(X\times 2^Y,T\times S)$ is distal, $\pi_1$ is a distal extension. By (3) the maximal equicontinuous factor of $(W,T\times S)$ consists of fixed points. Let $\pi_2: Z\ra 2^Y$ be the projection to the second coordinate,
i.e. $\pi_2(x,A)=A$.
Then $\Y_J=\pi_2(W)$ is equicontinuous as $2^Y$ is equicontinuous, which means that $J[x_0]$ is a fixed point of $(2^Y,S)$. Thus $J[x_0]=Y$, a contradiction with $J[x_0]\neq Y$. Thus $X\perp Y$.

\medskip
(1) $\Longrightarrow$ (4). Let $(X',T')$ be an almost one-to-one extension of $(X,T)$. Then $(X',T')\bot \mathcal{M}$ by Lemma \ref{disjoint for almost 1-1}. Hence, $(X',T')$ has DDMS-property by Lemma \ref{dense-mi}.

\medskip
(4) $\Longrightarrow$ (5) is clear.  (5)  $\Longrightarrow$ (3) This follows from Theorem \ref{ge-sufficient-condition}.
\end{proof}

\subsection{Proof of Theorem \ref{Main-thm-B}} Before proving Theorem \ref{Main-thm-B}, we need some lemmas. 
The first lemma was appeared first in \cite[Corollary 2.1]{HLSY03} and for completeness we give a proof here.
\begin{lem} \label{lift-rp}Let $\pi:(X,T)\ra (Y,S)$ be an extension between two systems. Then for any $(y,y')\in Q(Y)$ there is some $w\in Y$ such that  $(y, w), (w,y')\in Q(Y)$ and both $(y,w)$ and $(w,y')$ have $\pi$-lifts.
\end{lem}
\begin{proof}
Let $(y,y')\in Q(Y)$. Then there are  subsequences $\{y_i\}_{i\in \mathbb{N}}, \{y_i'\}_{i\in \mathbb{N}}$ in $Y$ and $\{n_i\}_{i\in \mathbb{N}}$ in $\mathbb{Z}$  such that $y_i\ra y, y_i'\ra y'$ and $S^{n_i}y_i\ra z, S^{n_i}y'_i\ra z$ for some $z\in Y$, as $i\rightarrow \infty$.

Without loss of generality, we may assume that $S^{-n_i}z\ra w$ for some $w\in Y$, as $i\rightarrow \infty$. Since
$\lim_{i\rightarrow \infty}\rho(S^{n_i}y_i, S^{n_i}S^{-n_i}z)=\lim_{i\rightarrow \infty}\rho(S^{n_i}y_i, z)=0$
we have $(y,w)\in Q(Y)$ and similarly we have $(w,y')\in Q(Y)$, where $\rho(\cdot,\cdot)$ is the metric on $Y$.

Now take $x_i\in \pi^{-1}(S^{n_i}y_i)$. Without loss of generality, assume  $x_i\ra x$  for some $x\in X$, and that $T^{-n_i}x_i\ra \hat{x}$  and $T^{-n_i}x\ra u$ for some $\hat{x},u\in X$, as $i\rightarrow \infty$. Then it is clear that $(\hat{x},u)\in Q(X)$.
Note that
$\pi(x)=\lim\limits_{i\rightarrow +\infty}\pi(x_i) =\lim \limits_{i\rightarrow +\infty} S^{n_i}(y_i)=z.$ Hence
$$\pi(u)=\lim_{i\rightarrow \infty}\pi(T^{-n_i}x)=\lim_{i\rightarrow \infty}S^{-n_i}z=w.$$
Since
$\lim_{i\rightarrow \infty}\pi(T^{-n_i}x_i)=\lim_{i\rightarrow \infty}y_i=y$,
one has $\pi(\hat{x})=y$. Thus we have shown that $(y,w)\in Q(Y)$ and it has a $\pi$-lift $(\hat{x},u)\in Q(X)$. Similarly, $(w,y')\in Q(Y)$ and it also has a $\pi$-lift.
\end{proof}

\begin{lem}\label{continuous of orb in semisimple}
Let $(X,T)$ be a semi-simple system.  If $\overline{\O} :X\ra 2^X$, $x\mapsto \overline{\O}(x,T)$ is continuous at point $x_0\in X$,  then $\overline{\O} $ is also continuous at every point in $\overline{\O}(x_0,T)$.
\end{lem}

\begin{proof}
Assume that $\overline{\O}$ is continuous at $x_0\in X$.
Let $y\in \overline{\O}(x_0,T)$. To show  $\overline{\O} $ is also continuous at  $y$, take $\{y_n\}$ in $X$ with $y_n\ra y$. Since $\overline{\O}$ is l.s.c. and $\overline{\O}(x_0,T)$ is minimal, for each $\epsilon>0$ one has
$\overline{\O}(x_0,T)=\overline{\O}(y,T)\subset B_{\epsilon}(\overline{\O}(y_n,T))$
when $n\in \mathbb{N}$ is large enough.
Thus we can find $y_n'\in \overline{\O}(y_n,T)$ for $n\in \mathbb{N}$ such that $\lim \limits_{n\rightarrow +\infty}y_n'=x$. Then the continuity of $\overline{\O}$ at $x_0$ implies that $\overline{\O}(x_0,T)= \lim \limits_{n\rightarrow +\infty} \overline{\O}(y'_n,T)$ in $2^X$. But $\overline{\O}(y'_n,T)= \overline{\O}(y_n,T)$ because $\overline{\O}(y_n,T)$ in the semi-simple system $X$ is  minimal. Thus we have
\[\overline{\O}(y,T)=\lim \limits_{n\rightarrow +\infty}\overline{\O}(y_n,T)\]
in $2^X$, and the continuity of $\overline{\O}$ at $y$ is followed.
\end{proof}

Let $Q^2(X)=\{(x,z)\in X\times X: \text{ there exists } y\in X \text{ such that }(x,y), (y,z)\in Q(X)\}$. Now we are going to prove Theorem \ref{Main-thm-B}.
\begin{proof}[Proof of Theorem \ref{Main-thm-B}] Let  $(X,T)$ be a distal system.
We show the following statements are equivalent:
\begin{enumerate}
\item $X\perp\mathcal{M}$.
\item For each minimal subsystem $W$ of $X$ we have $W\times W\subset Q(X)$.
\item There are minimal subsystems $\{X_i\}_{i\in\mathbb{N}}$ with $\bigcup_{i=1}^{\infty} X_i$ is dense in $X$, and for each $i\in\N$, $X_i\times X_i\subset Q(X)$.
\end{enumerate}

\medskip

(1) $\Longrightarrow$ (2). This follows from Corollary \ref{dis-sufficient-condition}.

(2) $\Longleftrightarrow$ (3). It is clear that (2) implies (3). Now assume that (3) holds. For any $x\in X$ there is a sequence $\{x_i\}_{i\in\N}\subset \cup_{j\in\N}M_j$ such that $\lim x_i=x$. Without loss of generality assume that $\overline{\O}(x_i,T)\ra M$. It is clear that $M\supset \overline{\O}(x,T)$. This implies that $M\times M\subset Q(X)$ and hence $\overline{\O}(x,T)\times \overline{\O}(x,T)\subset Q(X)$.

(2) $\Longrightarrow$ (1). Now we assume that (2) holds. Let $(Y,T)$ be an almost one-to-one distal extension of $(X,T)$ and $\pi:Y\ra X$  be the factor map.

Let $\Omega_1$ be the dense $G_\delta$ set $\{y\in Y: \pi^{-1}\pi(y)=\{y\}\}$ and $\Omega_2$ be the dense $G_\delta$ set of the continuous points
of $\overline{\O}: Y\rightarrow 2^Y, y\mapsto \overline{\O(y,T)}$. Set $\Omega=\Omega_1\cap \Omega_2$.

Now let $z\in \Omega$ and $y,y'\in \overline{\O}(z,T)$. We are going to show that $(y,y')\in Q^2(Y)$.
Since (2) holds, $(\pi(y)$, $\pi(y')\in Q(X)$. By Lemma \ref{lift-rp}, there is $x\in X$ with $(\pi(y),x),$ $(x,\pi(y'))\in Q(X)$
and both of them have $\pi$-lifts. So, by Lemma \ref{easylemmas} there are $y_1,y_2\in Y$ with $\pi(y_1)=\pi(y_2)=x$ with $(y,y_1),(y_2,y')\in Q(Y)$.

By Lemma \ref{continuous of orb in semisimple}, each point of $\overline{\O}(z,T)$ is a continuous point of the map $\overline{\O}$.
We conclude that $y_1,y_2\in \overline{\O}(z,T)$. To be precise, we show $y_1\in \overline{\O}(z,T)$ and $y_2\in \overline{\O}(z,T)$ follows by the same proof. Let $x_i\ra y_1, z_i\ra y$ and $n_i\in\Z$ with $\rho(T^{n_i}x_i,T^{n_i}z_i)\ra 0$. Let $\ep>0$. Then there is $\delta>0$ such that if $\rho(z_1,z_2)<\delta$ with $z_1\in \overline{\O}(z,T)$ then $\rho_H(\overline{\O}(z_1,T), \overline{\O}(z_2,T) )<\ep$. This implies that if $i$ is large, then $\rho_H(\overline{\O}(x_i,T), \overline{\O}(z,T))<\ep$. It follows that $y_1\in \overline{\O}(z,T)$, as $\ep$ is arbitrary.


Now Lemma \ref{easylemmas} implies that $y_1=y_2$ and thus $(y,y')\in Q^2(X)$. That is, $\overline{\O}(z,T)\times \overline{\O}(z,T)\subset Q^2(X)$. It is clear that $Q^2(X)\in \mathcal{A}(Q(X))$. This implies that the maximal equicontinuous factor of $(Y,T)$ is consisting of fixed points (by the density of $\Omega$) which in turn implies that $X\perp \mathcal{M}$ by Theorem \ref{equi-ll}(3).
\end{proof}

\section{Some open questions}\label{sec-open}

In the section, we  pose several questions. It has been showed in \cite{BGR23, GLR24} that the product property on disjointness holds in ergodic theory, that is if two measure preserving systems are both disjoint from every ergodic system then so is their product. In the appendix, we will give an alternative proof of this result. But this remains open in topological dynamics. Even under the transitivity condition, whether the product property holds for any two transitive systems remains open as asked in \cite{DSY12} and \cite{HSY20}. So we leave the following open question.
\begin{ques}
Let $(X,T)$ and $(Y,S)$ be two systems that are both disjoint from any minimal systems. Is the product $(X\times Y, T\times S)$ also disjoint from any minimal systems?
\end{ques}
We need to point out that in \cite{HSY20}, the authors show the following product property for transitive systems:
Let $(X,T)$ be a transitive system. If $(X,T)\perp \mathcal{M}$ then $(X^{(n)},T^{(n)})\perp \mathcal{M}$ for each $n\in\mathbb{N}$, where
$X^{(n)}=X\times \cdots \times X \,\text{($n$-times)}$ and $T^{(n)}=T\times \cdots\times T \,\text{($n$-times)}$.

As Theorem \ref{dense-mi} shows that if a system is disjoint from all minimal systems then it has dense minimal points. On the other hand, Theorem \ref{thm general} tells us that it has a topologically ergodic decomposition for an almost one-to-one extension. But we do not know whether there are many maximal transitive subsystem that has dense minimal points. We ask the following question.

\begin{ques}\label{many M-sys}
Suppose that $(X,T)$ is a system that is disjoint from all minimal systems. Is there a dense $G_{\delta}$ set $\Omega$ in $X$ such that for each $x\in \Omega$, the subsystem $\overline{\O}(x, T)$ has dense minimal points?
\end{ques}



\bigskip

\appendix

\section{Measure preserving systems disjoint from all ergodic automorphisms}\label{sec-appendix}

In this appendix, we will show that a measure preserving system is disjoint from every ergodic automorphism iff it is disjoint from almost every its ergodic component.

First we need some preparations.
Let $(X,\X,\mu,T)$ be an invertible probability measure preserving system. Let
$$\mu=\int_{\overline{X}}\mu_{\overline{x}}d\overline{\mu}(\overline{x})$$
be the ergodic decomposition of $T$. We can set $X=\overline{X}\times Z$, where both $\overline{X}, Z$ can be chosen as Polish spaces (see \cite[Subsection 1.3 ]{GLR24} or \cite[Section 2]{GW24}). Let $T_{\overline{x}}=T_{\{\overline{x}\}\times Z}$, or more precisely $T(\overline{x},z)=(\overline{x},T_{\overline x}z)$. Then  $(Z, \mu_{\overline{x}}, T_{\overline{x}})$ is the ergodic components of $T$. Let
$$D=\{(\overline{x}_1, \overline{x}_2)\in \overline{X}\times \overline{X}: T_{\overline{x}_1}\perp T_{\overline{x}_2}\}.$$
Then G\'{o}rska-Lema\'{n}czyk-de la Rue' theorem says that $(X,\X,\mu,T)$ is disjoint from every ergodic system iff $\overline{\mu}\times \overline{\mu}(D)=1$ (\cite[Theorem 3.1]{GLR24}).

To prove Theorem \ref{thm-A2}, we need the following lemma (compare it with Theorem \ref{hy-general}).

\begin{lem}[{\cite[Corollary 3.6, Lemma 3.7]{GLR24}\label{lem-GLR}}]
Let $(X,\X,\mu,T)$ be an invertible probability measure preserving system. Let
$\displaystyle \mu=\int_{\overline{X}}\mu_{\overline{x}}d\overline{\mu}(\overline{x})$
be the ergodic decomposition of $T$ described as above.
Let $(Y,\Y,\nu,S)$ be an invertible probability measure preserving system. Then the set $\{\overline{x}\in \overline{X}: T_{\overline{x}}\perp S\}$ is Borel. And if $\overline{\mu}\left(\{\overline{x}\in \overline{X}: T_{\overline{x}}\not\perp S\}\right)>0$, then $T\not \perp S$.
\end{lem}

Now we have the following description for the class ${\rm Erg}^\perp$ of automorphisms disjoint from all ergodic automorphisms..

\begin{thm}\label{thm-A2}
Let $(X,\X,\mu,T)$ be an invertible probability measure preserving system. Let
$\displaystyle \mu=\int_{\overline{X}}\mu_{\overline{x}}d\overline{\mu}(\overline{x})$
be the ergodic decomposition of $T$ described as above. Then $(X,\X,\mu,T)$ is disjoint from all ergodic automorphisms iff it is disjoint from almost all its ergodic components, i.e. $T\perp T_{\overline{x}}$ for $\overline{\mu}$-a.e. $\overline{x}\in \overline{X}$.
\end{thm}

\begin{proof}
If $(X,\X,\mu,T)$ is disjoint from every ergodic system, then it is clear that $T\perp T_{\overline{x}}$ for $\overline{\mu}$-a.e. $\overline{x}\in \overline{X}$ as $T_{\overline{x}}$ is ergodic.

Now we show the converse. Assume that  $T\perp T_{\overline{x}}$ for $\overline{\mu}$-a.e. $\overline{x}\in \overline{X}$. For each such $\overline{x}_1\in \overline{X}$, by Lemma \ref{lem-GLR}, the set $D_{\overline{x}_1}=\{\overline{x}_2\in \overline{X}: T_{\overline{x}_1}\perp T_{\overline{x}_2}\}$ is Borel, and $\overline{\mu}(D_{\overline{x}_1})=1$. That is, for $\overline{\mu}$-a.e. $\overline{x}_1\in \overline{X}$, $\overline{\mu}(D_{\overline{x}_1})=1$. By Fubini's Theorem, we have $\overline{\mu}\times \overline{\mu}(D)=1$.
Thus by G\'{o}rska-Lema\'{n}czyk-de la Rue' theorem, $(X,\X,\mu,T)$ is disjoint from every ergodic system. The proof is complete.
\end{proof}

Next we use G\'{o}rska-Lema\'{n}czyk-de la Rue' theorem to show that the class ${\rm Erg}^\perp$ is closed under products.

\begin{thm}[{\cite[Theorem 30]{BGR23},\cite[Theorem 6.1]{GW24}}]\label{thm-A3}
If both $(X,\X,\mu,T)$ and $(Y,\Y,\mu,S)$ are disjoint from all ergodic automorphisms, then so does their product system.
\end{thm}

\begin{proof}
Let
$\displaystyle \mu=\int_{\overline{X}}\mu_{\overline{x}}d\overline{\mu}(\overline{x})$
be the ergodic decomposition of $T$, and let
$\displaystyle \nu=\int_{\overline{Y}}\mu_{\overline{y}}d\overline{\nu}(\overline{y})$
be the ergodic decomposition of $S$. Let $\{T_{\overline{x}}\}_{\overline{x}\in \overline{X}}$ and $\{S_{\overline{y}}\}_{\overline{y}\in \overline{Y}}$ be the ergodic components of $T$ and $S$ respectively.

For $\overline{\mu}$-a.e. $\overline{y}\in \overline{Y}$, $S_{\overline{y}}$ is ergodic. By $T\perp S_{\overline{y}}$ and Lemma \ref{lem-GLR}, $\overline{\mu}\left(\{\overline{x}\in \overline{X}: T_{\overline{x}}\perp S_{\overline{y}}\}\right)=1$. Thus by Fubini's Theorem, for $\overline{\mu}\times \overline{\nu}$-a.e. $({\overline{x}},\overline{y})$, $T_{\overline{x}}\perp S_{\overline{y}}$, and therefore their product is ergodic. It follows that
$$\mu\times \nu=\int_{\overline{X}\times \overline{Y}}\mu_{\overline{x}}\times \nu_{\overline{y}} d(\overline{\mu}\times \overline{\nu})(\overline{x},\overline{y})$$
is the ergodic decomposition of $T\times S$.

Now for $\overline{\mu}\times \overline{\nu}$-a.e. $({\overline{x}_1},\overline{y}_1)$,
$T_{\overline{x}_1}\times S_{\overline{y}_1}$ is ergodic. And by $T\perp (T_{\overline{x}_1}\times S_{\overline{y}_1})$ and Lemma \ref{lem-GLR}, for $\overline{\mu}$-a.e. $\overline{x}_2\in \overline{X}$, we have that $T_{\overline{x}_2}\perp (T_{\overline{x}_1}\times S_{\overline{y}_1})$. Thus $(T_{\overline{x}_1}\times S_{\overline{y}_1})\times T_{\overline{x}_2}$ is ergodic. By $S\perp (T_{\overline{x}_1}\times S_{\overline{y}_1})\times T_{\overline{x}_2}$ and Lemma \ref{lem-GLR}, for $\overline{\nu}$-a.e. $\overline{y}_2\in \overline{Y}$, we have that $S_{\overline{y}_2}\perp (T_{\overline{x}_1}\times S_{\overline{y}_1})\times T_{\overline{x}_2}$, i.e. $(T_{\overline{x}_1}\times S_{\overline{y}_1})\perp (T_{\overline{x}_2}\times S_{\overline{y}_2})$. By Fubini's Theorem, we have that
For $\overline{\mu}\times \overline{\nu}$-a.e. $({\overline{x}_1},\overline{y}_1), ({\overline{x}_2},\overline{y}_2)\in \overline{X}\times \overline{Y}$,
$$(T_{\overline{x}_1}\times S_{\overline{y}_1})\perp (T_{\overline{x}_2}\times S_{\overline{y}_2}).$$
By  G\'{o}rska-Lema\'{n}czyk-de la Rue' theorem, $(X\times Y, \X\times \Y,\mu\times \nu, T\times S)$ is disjoint from all ergodic automorphisms. The proof is complete.
\end{proof}

\end{document}